\documentclass[11pt]{article}
\usepackage{graphicx,wrapfig,subfigure}
\usepackage{amsmath}
\usepackage{amsthm}
 \usepackage{amssymb}
\usepackage{setspace}
\usepackage{leftidx}
\usepackage{graphicx}
\usepackage{amsfonts}
\usepackage{color, epstopdf}
\usepackage{cite}
\usepackage{array,setspace,diagbox}
\usepackage{booktabs}
\usepackage [latin1]{inputenc}
\usepackage{multirow}
\usepackage{enumitem}
\usepackage{algpseudocode,algorithm,algorithmicx}%
\usepackage{appendix}%
\usepackage{upgreek}
\usepackage{enumerate}
\usepackage{enumitem}

\newtheorem{theorem}{Theorem}[section]
\newtheorem{lemma}{Lemma}[section]

\newtheorem{example}{Example}
\newtheorem{remark}{Remark}

\newcommand{\diff}{\triangledown_{\!\tau}}

\newcommand{\defeq}{:=}
\newcommand{\zd}{\,\mathrm{d}}
\newcommand{\abs}[1]{\left|#1\right|}

\newcommand{\bra}[1]{\left(#1\right)}
\newcommand{\brab}[1]{\big(#1\big)}
\newcommand{\braB}[1]{\Big(#1\Big)}
\newcommand{\brat}[1]{(#1)}
\newcommand{\kbra}[1]{\left[#1\right]}
\newcommand{\kbrab}[1]{\big[#1\big]}
\newcommand{\kbraB}[1]{\Big[#1\Big]}
\newcommand{\kbrat}[1]{[#1]}
\newcommand{\myinner}[1]{\left\langle#1\right\rangle}
\newcommand{\myinnerb}[1]{\big\langle#1\big\rangle}

\newcommand{\mynorm}[1]{\left\|#1\right\|}
\newcommand{\mynormb}[1]{\big\|#1\big\|}
\newcommand{\mynormB}[1]{\Big\|#1\Big\|}

\newcommand{\Rcc}[1]{{}^{r}\!{#1}}

\begin{document}
\title{Discrete gradient structure of a second-order variable-step method
 for nonlinear integro-differential models}
\author{
	Hong-lin Liao\thanks{ORCID 0000-0003-0777-6832. College of Mathematics,
		Nanjing University of Aeronautics and Astronautics,
		Nanjing 211106, China; Key Laboratory of Mathematical Modeling
		and High Performance Computing of Air Vehicles (NUAA), MIIT, Nanjing 211106, China.
		Emails: liaohl@nuaa.edu.cn and liaohl@csrc.ac.cn.
		This author's work is supported by NSF of China
		under grant number 12071216.}
	\and
	Nan Liu\thanks{College of Mathematics, Nanjing University of Aeronautics and Astronautics,
		Nanjing 211106, China. Email: liunan@nuaa.edu.cn.}	
	\and
	Pin Lyu\thanks{ORCID 0000-0002-6947-6838. School of Mathematics, Southwestern University of
        Finance and Economics, Chengdu 611130, China. Email: plyu@swufe.edu.cn.
		This author is supported by NSF of China under grant number 12101510,
        NSF of Sichuan Province under grant number 2022NSFSC1789, and Guanghua Talent Project of SWUFE.}	
	}
\date{\today}
\maketitle
\normalsize

\begin{abstract}
The discrete gradient structure and the positive definiteness of
discrete fractional integrals or derivatives are fundamental
to the numerical stability in long-time simulation of nonlinear integro-differential models.
We build up a discrete gradient structure for a class of second-order variable-step approximations
of fractional Riemann-Liouville integral and fractional Caputo derivative.
Then certain variational energy dissipation laws at discrete levels of the
corresponding variable-step Crank-Nicolson type methods
are established for time-fractional Allen-Cahn
and time-fractional  Klein-Gordon type models.
They are shown to be asymptotically compatible with the associated energy laws
of the classical Allen-Cahn and Klein-Gordon equations in the associated fractional order limits.
Numerical examples together with an adaptive time-stepping procedure
 are provided to demonstrate the effectiveness of our second-order methods.
  \noindent{\bf AMS subject classiffications.}\;\; 35Q99, 65M06, 65M12, 74A50
\end{abstract}

\section{Introduction}
\setcounter{equation}{0}

In the past few decades, linear and nonlinear integro-differential equations
attract great interests in a wide range of disciplines in science and engineering \cite{Brunner2004book,GolmankhanehGolmankhanehBaleanu:2011,Mainardi2010,MetzlerKlafter:2000,Nigmatullin-pb1984}.
Typically, diffusion equations with fractional derivatives
and fractional integrals have become widely-used models describing anomalous
diffusion processes where the mean squared displacement
scales as a fractional power of time. These models
exhibit multi-scaling time behaviour, which makes them suitable for the description of
different diffusive regimes and characteristic crossover dynamics in complex systems
\cite{Mainardi2010,MetzlerKlafter:2000,Nigmatullin-pb1984}.
They are always formulated in the integral form, including the Riemann-Liouville fractional integral
\begin{align}\label{def: RL integral}
	(\mathcal{I}_{t}^{\beta}w)(t)
	:=\int_{0}^{t}\omega_{\beta}(t-s)w(s)\zd{s}\quad\text{with $\omega_{\beta}(t):=t^{\beta-1}/\Gamma(\beta)$ for $\beta>0$,}
\end{align}
and the fractional Caputo derivative
\begin{align}\label{def:Caputo derivative}
	(\partial_{t}^{\alpha}w)(t)
	:=(\mathcal{I}_{t}^{1-\alpha}w^\prime)(t)
	=\int_{0}^{t}\omega_{1-\alpha}(t-s)w'(s)\zd{s}\quad\text{for  $0<\alpha<1$.}
\end{align}

In capturing the multi-scale behaviors in
many of integro-differential equations, such as the time-fractional phase field models \cite{JiLiaoGongZhang:2020,JiZhuLiao:2022,LiaoTangZhou:2020,LiaoTangZhou:2021,QuanTangYang-2020csiam,
	QuanTangYang-2020,QuanTangWangYang:2022,TangYuZhou:2019}
and nonlinear fractional wave models \cite{AlsaediAhmadKirane:2015,AdolfssonEnelundLarsson:2003,
	CuestaLubichPalencia:2006,CuestaPalencia:2003,GolmankhanehGolmankhanehBaleanu:2011,
	LubichSloanThomee:1996,LyuVong:2022JSC,MustaphaMustapha:2010},
adaptive time-stepping strategies, namely, small time steps are utilized when the solution varies
rapidly and large time steps are employed otherwise, are practically useful
\cite{Brunner2004book,LiaoTangZhou:2020,LiaoZhuWang:2022NMTMA,LyuVong:2022JSC,McLeanMustapha:2007,
	McLeanThomee:1996,MustaphaMustapha:2010,Mustapha-sinum2020,Mustapha McLean:2013,MustaphaSchotzau:2014}.
It requires practically and theoretically reliable (stable and convergent)
time-stepping methods on general setting of time step-size variations
\cite{JiLiaoGongZhang:2020,JiZhuLiao:2022,LiaoTangZhou:2020,LiaoTangZhou:2021, LiaoZhuWang:2022NMTMA}.


This work is concerned with a class of second-order approximations
with unequal time-steps of fractional Riemann-Liouville integral and fractional Caputo derivative.
For a given time $T>0$ and a positive integer $N$,
consider the time levels $0=t_{0}<t_{1}<\cdots<t_{k}<\cdots<t_{N}=T$
with the step sizes $\tau_{k}:=t_{k}-t_{k-1}$ for $1\le k\le N$.
The maximum step size is denoted by $\tau:=\max_{1\le k\le N}\tau_k$
and the local time-step ratio $r_k:=\tau_k/\tau_{k-1}$ for $2\le k\le N$.
Given a grid function $w^k=w(t_k)$, define
\begin{align*}
	w^{k-\frac12}:=(w^k+w^{k-1})/2,\quad \diff w^k:=w^k-w^{k-1},\quad \partial_\tau w^k:=\diff w^k/\tau_k\quad
	\text{for $k\geq1$}.
\end{align*}
Let $(\Pi_{0,k}w)(t)$ be the constant interpolant of a function $w$ at two nodes $t_{k-1}$ and $t_{k}$,
and define the piecewise approximation $\Pi_{0}w:=\Pi_{0,k}w$  so that $(\Pi_{0}w)(t)=w^{k-\frac12}$
for $t_{k-1}<{t}\leq t_{k}$ and $k\geq1$. The integral averaged formula
(also called Crank-Nicolson approximation in \cite{McLeanThomee:1993,McLeanThomee:1996,McLeanMustapha:2007})
of fractional Riemann-Liouville integral \eqref{def: RL integral} reads
\begin{align}\label{def: IAF formula}
	(\mathcal{I}_{\tau}^{\beta}w)^{n-\frac12}
	:=\frac{1}{\tau_n}\int_{t_{n-1}}^{t_n}\int_0^{t}
	\omega_{\beta}(t-s)(\Pi_{0}w)\zd{s}\zd{t}
	\triangleq\sum_{k=1}^{n}a_{n-k}^{(\beta,n)}\tau_kw^{k-\frac12},
\end{align}
where the associated discrete  kernels $a_{n-k}^{(\beta,n)}$ are defined by
\begin{align}\label{def: IAF formula kernel}
	a_{n-k}^{(\beta,n)}&:=\frac{1}{\tau_n\tau_{k}}\int_{t_{n-1}}^{t_n}\int_{t_{k-1}}^{\min\{t,t_k\}}
	\omega_{\beta}(t-s)\zd{s}\zd{t}\quad\text{for $1\leq{k}\leq{n}$}.
\end{align}
Let $(\Pi_{1,k}w)(t)$ denote the linear interpolant at two nodes $t_{k-1}$ and $t_{k}$,
and define the piecewise approximation $\Pi_{1}w:=\Pi_{1,k}w$ so that  $(\Pi_{1}w)'(t)=\partial_{\tau}w^{k}$
for $t_{k-1}<{t}\leq t_{k}$ and $k\geq1$. The integral averaged  formula
(also called L1$^+$ formula \cite{JiLiaoGongZhang:2020}) of fractional Caputo derivative \eqref{def:Caputo derivative} is
\begin{align}\label{def: L1+ formula}
	(\partial_{\tau}^{\alpha}w)^{n-\frac12}
	:=\frac{1}{\tau_n}\int_{t_{n-1}}^{t_n}\int_0^{t}
	\omega_{1-\alpha}(t-s)(\Pi_{1}w)^\prime(s)\zd{s}\zd{t}
	\triangleq\sum_{k=1}^{n}a_{n-k}^{(1-\alpha,n)}\diff w^{k},
\end{align}
where the associated discrete kernels $a_{n-k}^{(1-\alpha,n)}$ are defined by \eqref{def: IAF formula kernel}, namely,
\begin{align}\label{def: L1+ formula kernel}
	a_{n-k}^{(1-\alpha,n)}=\frac{1}{\tau_n\tau_{k}}\int_{t_{n-1}}^{t_n}\int_{t_{k-1}}^{\min\{t,t_k\}}
	\omega_{1-\alpha}(t-s)\zd{s}\zd{t}\quad\text{for $1\leq{k}\leq{n}$}.
\end{align}
Both of these approximations come from
the so-called positive-semidefinite-preserving approach
for approximating convolution integrals, in which the numerical approximations are 
designed such that
the corresponding real quadratic form is a discrete analogue to
the non-negative definiteness of continuous kernels, see \cite{McLeanThomee:1993,McLeanThomee:1996,McLeanMustapha:2007,QuanTangYang-2020csiam},
\begin{align}\label{ineq: continuous semi-positive}
	2\,\mathcal{I}_t^{1}\brab{w\,\mathcal{I}_t^{\beta} w}(t)
	=&2\int_0^tw(\mu)\zd\mu\int_0^{\mu}\omega_{\beta}(\mu-s)w(s)\zd s\nonumber\\
	=&\int_0^t\int_0^tw(s)w(\mu)\omega_{\beta}(\abs{\mu-s})\zd\mu\zd s\geq0
	\quad\text{for $t>0$ and $w\in C[0,T]$.}
\end{align}
Recently, the regularity condition $w\in C[0,T]$ was updated by 
Tang et al. \cite[Lemma 2.1 and Corollary 2.1]{TangYuZhou:2019}, that is, 
the semipositive definiteness \eqref{ineq: continuous semi-positive} holds for $w \in L^p(0,T)$ with $p\ge \frac{2}{1+\beta}$ for $0<\beta<1$, which permits some weakly singular functions like $w=O(t^{\beta-1})$ such that the discrete kernels of L1$^+$ formula \eqref{def: L1+ formula} can 
naturally preserve the non-negative definiteness.
As seen, the discrete integral \eqref{def: IAF formula} and
the discrete derivative \eqref{def: L1+ formula} are different
in the regularity requirement of $w$ to ensure the non-negative definiteness.
However, the same expression of \eqref{def: IAF formula kernel} and \eqref{def: L1+ formula kernel}
urges us to explore whether we can determine the positive definiteness of these discrete convolution kernels
without using the non-negative definiteness of continuous kernels.

Under a mild step-ratio condition, we establish the following \emph{discrete gradient structure} (DGS)
for the discrete convolution kernels $a_{n-k}^{(\beta,n)}$ defined by
in \eqref{def: IAF formula kernel}, that is,
\begin{align}\label{eq: DGS equality}
	2w_n\sum_{j=1}^na^{(\beta,n)}_{n-j}w_j=&\,
	\sum_{j=1}^{n}\mathsf{p}_{n-j}^{(\beta,n)}v_j^2
	-\sum_{j=1}^{n-1}\mathsf{p}_{n-1-j}^{(\beta,n-1)}v_j^2\nonumber\\
	&\,+\sum_{j=1}^{n-1}\braB{\frac{1}{\Rcc{\mathsf{p}}_{n-j}^{(\beta,n)}}-\frac{1}{\Rcc{\mathsf{p}}_{n-j-1}^{(\beta,n)}}}
	\braB{\sum_{k=1}^{j}\Rcc{\mathsf{p}}_{n-k}^{(\beta,n)}(v_k-v_{k-1})}^2\quad\text{for $n\ge1$,}
\end{align}
where the sequence $v_j:=\sum_{\ell=1}^{j}\mathsf{a}_{j-\ell}^{(\beta,j)}w_{\ell}$ with
the modified kernels
\begin{align}\label{def: modified kernels}
	\mathsf{a}_{0}^{(\beta,n)}:=2 a_{0}^{(\beta,n)}\quad\text{and}\quad
	\mathsf{a}_{n-j}^{(\beta,n)}:=a_{n-j}^{(\beta,n)}\quad\text{for $1\le j\le n-1$.}
\end{align}
Here, $\mathsf{p}_{n-j}^{(\beta,n)}>0$ and $\Rcc{\mathsf{p}}_{n-j}^{(\beta,n)}>0$ are the associated
discrete (left-)complementary convolution (DCC) and right-complementary convolution (RCC)
kernels with respect to the modified kernels.  Note that this equality \eqref{eq: DGS equality}
provides a solution to the remaining question in \cite[Remark 4.1]{LiaoTangZhou:2020positive}.
Actually, it shows that the discrete convolution kernels in 
\eqref{def: IAF formula kernel} and \eqref{def: L1+ formula kernel}
are positive definite in the sense that
\begin{align}\label{ieq: positive definite}
	2\sum_{k=1}^nw_k\sum_{j=1}^ka^{(\beta,k)}_{k-j}w_j
	\ge \sum_{k=1}^{n}\mathsf{p}_{n-k}^{(\beta,n)}v_k^2>0\quad\text{for nonzero $\{w_{k}\}_{k=1}^n$.}
\end{align}

Another aim of this paper is to show that the resulting second-order Crank-Nicolson method
using the discrete integral \eqref{def: IAF formula} or
the discrete derivative \eqref{def: L1+ formula} inherits certain (variational) energy dissipation laws
at discrete time levels for time-fractional Allen-Cahn
and time-fractional Klein-Gordon models.
These energy dissipation laws are shown to be asymptotically compatible with the associated energy
dissipation (or conservation) laws
of the classical Allen-Cahn  and Klein-Gordon equations in the associated fractional order limits.

The rest of this paper is organized as follows.
In section 2, we derive the DGS of the discrete kernels \eqref{def: IAF formula kernel}.
The discrete energy dissipation law of the L1$^+$ scheme
for time-fractional Allen-Cahn model is addressed in section 3.
Section 4 presents a novel energy dissipation law of the time-fractional Klein-Gordon
equation and establishes the discrete counterpart.
Numerical examples are presented in section 5 to show the effectiveness of our time-stepping methods.

\section{Discrete gradient structure}
\setcounter{equation}{0}

\subsection{DGS for general kernels}

This section builds up the DGS \eqref{eq: DGS equality}.
 At first, we prove a continuous counterpart, 
where we use the Riemann-Liouville fractional derivative
${}^R\!\partial_t^{\beta}$ defined by
\begin{align*}
	{}^R\!\partial_t^{\beta} v\defeq \partial_t\mathcal{I}_t^{1-\beta}v
	\quad\text{for $0<\beta<1$}.
\end{align*}

\begin{lemma}\label{lem: countinous DGS equality}
	For $\beta\in(0,1)$ and an absolutely continuous function $w$, it holds that
	\begin{align*}
		2w(t)(\mathcal{I}_t^{\beta}w)(t)
		=&\,\brab{{}^R\!\partial_t^{\beta}v^2}(t)
		+\int_0^t\frac{\partial}{\partial\xi}\bra{\frac{1}{\omega_{1-\beta}(t-\xi)}}
		\bra{\int_0^{\xi}\omega_{1-\beta}(t-s)v'(s)\zd{s}}^2\zd \xi,
	\end{align*}
where $v=\mathcal{I}_t^{\beta}w$ and  $\frac{\partial}{\partial\xi}\brab{\frac{1}{\omega_{1-\beta}(t-\xi)}}=\Gamma(1-\beta)\Gamma(1+\beta)\omega_{\beta}(t-\xi)>0$, such that
\begin{align*}
	2\,\mathcal{I}_t\brab{w\mathcal{I}_t^{\beta}w}(t)\ge
	\brab{\mathcal{I}_t^{1-\beta}v^2}(t)>0\quad\text{for $v\not\equiv 0$.}
\end{align*}
\end{lemma}
\begin{proof} By the semigroup property, we have
	$w={}^R\!\partial_t^{\beta}v=\partial_t\mathcal{I}_t^{1-\beta}v$ and
	\begin{align*}
		w(t)(\mathcal{I}_t^{\beta}w)(t)=v(t)({}^R\!\partial_t^{\beta}v)(t).
	\end{align*}
	Since $v(0)=0$, one has ${}^R\!\partial_t^{\beta} v=\partial_{t}^{\beta}v$
	and ${}^R\!\partial_t^{\beta} v^2=\partial_{t}^{\beta}v^2$.	
	Consider the difference
	\begin{align*}
		J[v]:=&\,2v(t)({}^R\!\partial_t^{\beta}v)(t)
		-\brab{{}^R\!\partial_t^{\beta}v^2}(t)\\
		=&\,2v(t)\frac{\partial}{\partial t}\int_0^t\omega_{1-\beta}(t-s)v(s)\zd{s}
		-\frac{\partial}{\partial t}	\int_0^t\omega_{1-\beta}(t-s)v^2(s)\zd{s} \\
		=&\,2\int_0^t\omega_{1-\beta}(t-s) v'(s)\kbrab{v(t)-v(s)}\zd{s}
		=2\int_0^t\omega_{1-\beta}(t-s) v'(s)\int_s^tv'(\xi)\zd{\xi}\zd{s}\\
		=&\,2\int_0^tv'(\xi)\zd{\xi}
		\int_0^{\xi}\omega_{1-\beta}(t-s) v'(s)\zd{s},	
	\end{align*}
	where the integration order was exchanged in the last equality.
	By taking $$u(\xi):=\int_0^{\xi}\omega_{1-\beta}(t-s)v'(s)\zd{s}$$ with $u(0)=0$ and
	$u'(\xi)=\omega_{1-\beta}(t-\xi)v'(\xi)$, it is not difficult to derive that
	\begin{align}\label{Cont: integration by parts}
		J[v]=&\,2\int_0^t\frac{u'(\xi)u(\xi)}{\omega_{1-\beta}(t-\xi)}\zd{\xi}
		=\int_0^t\frac{\zd u^2(\xi)}{\omega_{1-\beta}(t-\xi)}
		=\int_0^t\frac{\partial}{\partial\xi}\bra{\frac{1}{\omega_{1-\beta}(t-\xi)}}u^2(\xi)\zd\xi.
	\end{align}
It leads to the claimed result and completes the proof.	
\end{proof}


To seek the discrete counterpart of  Lemma \ref{lem: countinous DGS equality},
we introduce some discrete tools for any discrete kernels $\{a_{n-j}^{(n)}\}_{j=1}^n$.
The discrete orthogonality convolution (DOC)
kernels $ {\theta_{n-k}^{(n)}} $ are defined by
\begin{align}\label{def: general DOC kernels}
	{\theta_{0}^{(n)}}:=\frac{1}{a_{0}^{(n)}}
	\quad \text{and}\quad
	{\theta_{n-k}^{(n)}}:=-\frac{1}{a_{0}^{(k)}}
	\sum_{j=k+1}^n {\theta_{n-j}^{(n)}}a_{j-k}^{(j)}
	\quad \text{for $1\le k\le n-1$.}
\end{align}
It is easy to check the following discrete orthogonality identities
\begin{align}\label{def: DOC orthogonal identity}
	\sum_{j=k}^n {\theta_{n-j}^{(n)}}a_{j-k}^{(j)}
	\equiv \delta_{nk}\quad\text{and}\quad
	\sum_{j=k}^n {a_{n-j}^{(n)}}\theta_{j-k}^{(j)}
	\equiv \delta_{nk}\quad\text{for $1\le k\le n$,}
\end{align}
where $\delta_{nk}$ is the Kronecker delta symbol with $\delta_{nk}=0$ if $k\neq n$.
With the DOC kernels, we define the discrete (left-)complementary
convolution (DCC) kernels
\begin{align}\label{def: general DCC kernels}
	p_{n-k}^{(n)}:=\sum_{j=k}^n {\theta_{j-k}^{(j)}}
	\quad \text{for $1\le k\le n$},
\end{align}
and the right-complementary convolution (RCC) kernels
\begin{align}\label{def: general RCC kernels}
	\Rcc{p}_{n-k}^{(n)}:=\sum_{j=k}^{n}{\theta_{n-j}^{(n)}}\quad\text{for $1\le k\le n$.}
\end{align}
By using the first orthogonality identity in \eqref{def: DOC orthogonal identity}, one can check that
\begin{align*}
	p_{0}^{(n)}a_{0}^{(n)}=1\quad\text{and}\quad
	\sum_{j=k}^{n}p_{n-j}^{(n)}a_{j-k}^{(j)}=\sum_{j=k}^{n-1}p_{n-j-1}^{(n-1)}a_{j-k}^{(j)}\quad\text{for $1\le k\le n-1$,}
\end{align*}
such that the DCC kernels $p_{n-j}^{(n)}$ are complementary with respect to  $a_{n-j}^{(n)}$ in the sense that
\begin{align}\label{def: DCC identity}
	\sum_{j=k}^{n}p_{n-j}^{(n)}a_{j-k}^{(j)}\equiv1	\quad\text{for $1\le k\le n$.}
\end{align}
In a similar fashion, by using the second orthogonality identity in \eqref{def: DOC orthogonal identity},
one can check that
\begin{align*}
	a_{0}^{(n)}\,\Rcc{p}_{0}^{(n)}=1\quad\text{and}\quad
	\sum_{j=k}^{n}a_{n-j}^{(n)}\Rcc{p}_{j-k}^{(j)}=\sum_{j=k+1}^{n}a_{n-j}^{(n)}\Rcc{p}_{j-k-1}^{(j)}\quad\text{for $1\le k\le n-1$.}
\end{align*}
So the kernels $a_{n-j}^{(n)}$ are complementary with respect to the RCC kernels $\Rcc{p}_{n-j}^{(n)}$ in the sense that
\begin{align}\label{def: RCC identity}
	\sum_{j=k}^{n}a_{n-j}^{(n)}\Rcc{p}_{j-k}^{(j)}\equiv1	\quad\text{for $1\le k\le n$.}
\end{align}

\begin{figure}[htb!]
	\centering
	\includegraphics[width=3.4in]{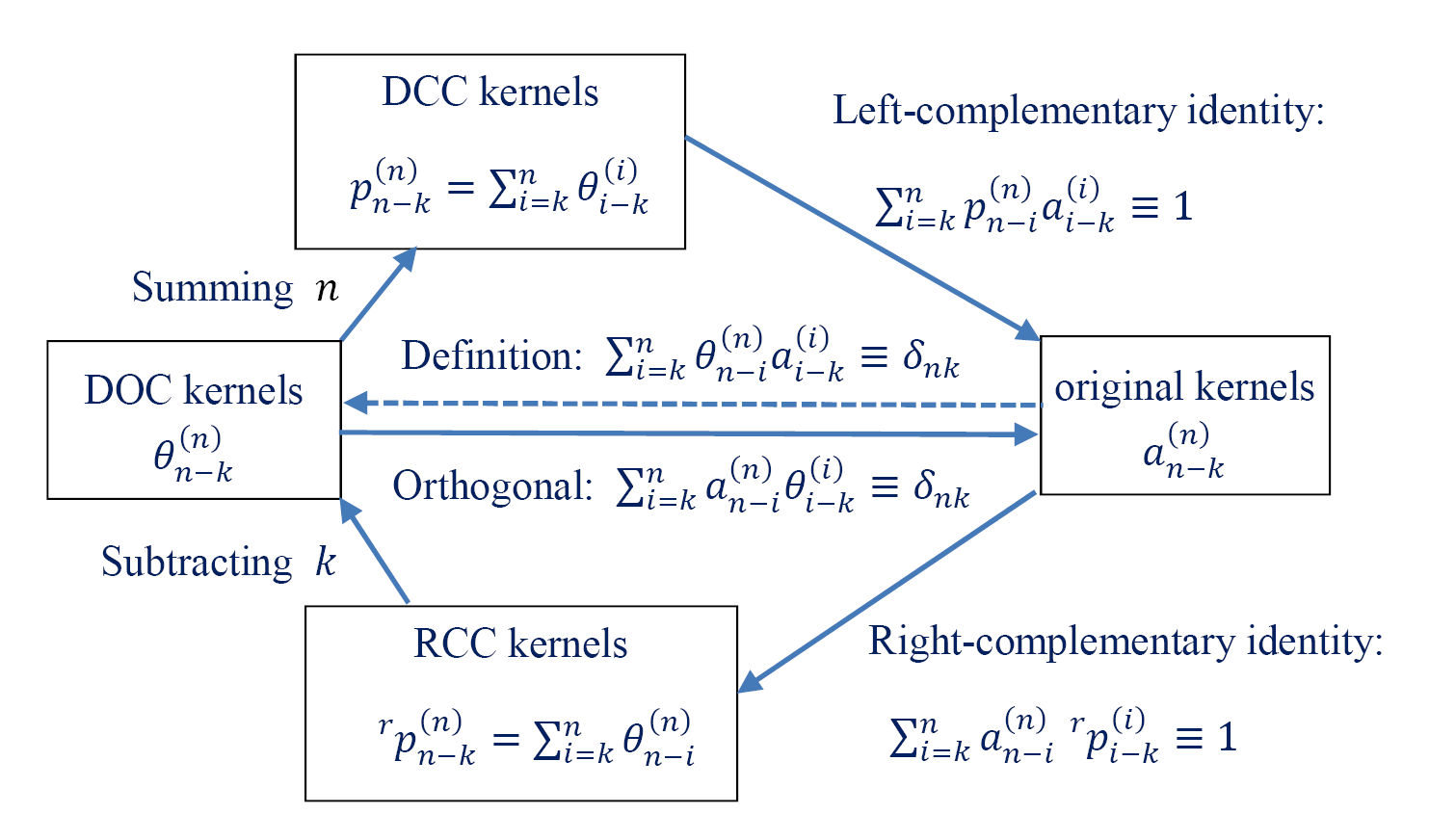}
	\caption{The relationship diagram of DOC, DCC and RCC kernels.}
	\label{fig: DOC, RCC and DCC relation}
\end{figure}

The above definitions and their connections are summarized in Figure \ref{fig: DOC, RCC and DCC relation}.
We note that the DOC kernels were originally introduced in \cite{LiaoZhang:2021}
for the analysis of variable-step BDF2 method of the first time derivative
and the DCC kernels were introduced in \cite{LiaoMcLeanZhang:2019}
to develop the discrete fractional Gr\"{o}nwall inequality for nonlinear subdiffusion equations with the fractional Caputo derivative. To the best of our knowledge, the RCC kernels
\eqref{def: general RCC kernels} are introduced here for the first time.
If	the given kernels $a_{n-j}^{(n)}$ simulate the continuous kernel $\omega_\beta$ 
	of the Riemann-Liouville integral \eqref{def: RL integral}, the above RCC and DCC kernels 
	are the discrete counterparts of the complementary kernel $\omega_{1-\beta}$ 
	in the sense that the complementary identities 
	\eqref{def: DCC identity} and \eqref{def: RCC identity} simulate the complementary (semigroup) property $\omega_{1-\beta}*\omega_\beta=1$ and $\omega_\beta*\omega_{1-\beta}=1$, respectively.
As is seen in the next two lemmas, the discrete properties of the DCC and RCC kernels are quite different.

\begin{lemma}\label{lem: DCC property}
	If the positive kernels $a^{(n)}_{j}$ are monotonically decreasing
	with respect to the subscript index $j$,
	that is,  $a^{(n)}_{j-1}>a_{j}^{(n)}$ for $1\le j\le n-1$, then the DCC kernels $p^{(n)}_{n-k}\ge0$.
\end{lemma}
\begin{proof}This result is obvious from the complementary identity \eqref{def: DCC identity}.
\end{proof}

\begin{lemma}\label{lem: RCC property}
	If the positive kernels $a^{(n)}_{j}$ are monotonically decreasing with respect to the superscript index $n$,
	$a^{(n-1)}_{j-1}>a_{j}^{(n)}$ for $1\le j\le n-1$,
	and satisfy a class of geometric-like convexity,
	$a_{j-1}^{(n-1)}a_{j+1}^{(n)}\ge a_{j}^{(n-1)}a_{j}^{(n)}$ for $1\le j\le n-2$,
	then the RCC kernels $\Rcc{p}_{j}^{(n)}$ in \eqref{def: general RCC kernels} are positive and monotonically decreasing
	with respect to $j$.
\end{lemma}
\begin{proof} The assumptions meet the condition of \cite[Lemma 2.3]{LiaoTangZhou:2020positive}, which
gives the following results
	on the associated DOC kernels $ {\theta_{n-k}^{(n)}}$,
	\begin{align*}
		\theta_{0}^{(n)}>0\quad\text{and}\quad \theta_{n-k}^{(n)}\le 0\;\;\text{ for $1\le k\le n-1$};\quad
		\text{but}\quad \sum_{k=1}^n {\theta_{n-k}^{(n)}}>0.
	\end{align*}
	Then the  definition \eqref{def: general RCC kernels} leads to the claimed result directly.
\end{proof}

The following theorem generalizes the result of \cite[Theorem 2.1]{JiZhuLiao:2022},
where the latter is valid for the well-known L1 formula of fractional Caputo derivative \eqref{def:Caputo derivative}.
\begin{theorem}\label{thm: DGS-general kernels}
	Let $n \ge 2$ and $\{\chi_{n-j}^{(n)}\}_{j=1}^n$ be a class of discrete convolution kernels.
	Consider the following auxiliary kernels for a constant $\sigma_{\min}\in[0,2)$,
	\[a_{0}^{(n)}:=(2-\sigma_{\min})\chi_{0}^{(n)}\quad\text{and}
	\quad a_{n-j}^{(n)}:=\chi_{n-j}^{(n)}\quad\text{for $1\le j\le n-1$.}\]
	Assume that the auxiliary kernels $a_{n-j}^{(n)}$ satisfy the following assumptions:
	\begin{description}	
		\item[(Row decrease)] \quad
		$\displaystyle a^{(n)}_{j-1}\ge a^{(n)}_{j}>0$\; for $1\le j\le n-1$;	
		\item[(Column decrease)] \quad
		$\displaystyle a^{(n-1)}_{j-1}>a_{j}^{(n)}$\; for $1\le j\le n-1$;		
		\item[(Logarithmic convexity)] \quad
		$\displaystyle a_{j-1}^{(n-1)}a_{j+1}^{(n)}\ge a_{j}^{(n-1)}a_{j}^{(n)}$\; for $1\le j\le n-2$.
	\end{description}
	Let $p_{n-j}^{(n)}$ and $\Rcc{p}_{n-j}^{(n)}$ be the associated DCC and RCC kernels, respectively,
	with respect to the modified kernels $a_{n-j}^{(n)}$.
	Then for any real sequence $\{w_k\}_{k=1}^n$, the following DGS holds,
	\begin{align*}
		2w_n\sum_{j=1}^n\chi_{n-j}^{(n)}w_j
		=&\,\sum_{k=1}^{n}p_{n-k}^{(n)}v_k^2-\sum_{k=1}^{n-1}p_{n-k-1}^{(n-1)}v_k^2
		+\sigma_{\min}\chi_{0}^{(n)}w_n^2\nonumber\\
		&+\sum_{j=1}^{n-1}\braB{\frac{1}{\Rcc{p}_{n-j}^{(n)}}-\frac{1}{\Rcc{p}_{n-j-1}^{(n)}}}
		\kbra{\sum_{k=1}^{j}\Rcc{p}_{n-k}^{(n)}(v_k-v_{k-1})}^2\quad \text{for $n\ge1$,}
	\end{align*}
	where $v_k:=\sum_{\ell=1}^{k}a_{k-\ell}^{(k)}w_{\ell}$ so that
	the convolution kernels $\chi_{n-k}^{(n)}$ are positive definite,
	\begin{align*}
		2\sum_{k=1}^nw_k\sum_{j=1}^k\chi_{k-j}^{(k)}w_j\geq
		\sum_{j=1}^{n}p_{n-j}^{(n)}v_j^2
		+\sigma_{\min}\sum_{k=1}^{n}\chi_{0}^{(k)}w_k^2\quad \text{for $n\ge1$}.
	\end{align*}
\end{theorem}

\begin{proof}Obviously, the auxiliary kernels $a_{n-j}^{(n)}$ fulfill the assumptions of
	Lemmas \ref{lem: DCC property} and \ref{lem: RCC property} so that the DCC kernels $p^{(n)}_{j}\ge0$
	and the RCC kernels $\Rcc{p}_{j}^{(n)}$ are positive and monotonically decreasing
	with respect to the subscript index $j$.

	For any real sequence $\{w_k\}_{k=1}^n$, let $v_0:=0$ and $v_j:=\sum_{k=1}^ja_{j-k}^{(j)}w_k$
	for $1\le j\le n$.
	With the help of the first discrete orthogonality identity in \eqref{def: DOC orthogonal identity}, 
	it is easy to derive that
	\begin{align*}
		w_k=\sum_{j=1}^k {\theta_{k-j}^{(k)}}v_j\quad\text{for $1\le k\le n$.}
	\end{align*}
	Then one applies the definitions \eqref{def: general DCC kernels} 
	and \eqref{def: general RCC kernels} to find
	\begin{align}\label{def: wn-vk}
		w_n=&\,\sum_{k=1}^{n}{\theta_{n-k}^{(n)}}v_k
		=\Rcc{p}_{0}^{(n)}v_n+\sum_{k=1}^{n-1}(\Rcc{p}_{n-k}^{(n)}-\Rcc{p}_{n-k-1}^{(n)})v_k
		=\sum_{k=1}^{n}\Rcc{p}_{n-k}^{(n)}\diff v_k,
	\end{align}
and
\begin{align}\label{def: vk-vk2}
	\sum_{k=1}^{n}\Rcc{p}_{n-k}^{(n)}\diff v_k^2=\sum_{k=1}^{n}{\theta_{n-k}^{(n)}}v_k^2 =\sum_{k=1}^{n}p_{n-k}^{(n)}v_k^2-\sum_{k=1}^{n-1}p_{n-k-1}^{(n-1)}v_k^2.
\end{align}

Fix~$n$ and consider the difference
\begin{align}\label{lemp: difference term Jn0}
	J_n:=&\,2w_n\sum_{j=1}^na_{n-j}^{(n)}w_j-\sum_{k=1}^{n}\Rcc{p}_{n-k}^{(n)}\diff v_k^2\\
	=&\,2v_n\sum_{k=1}^n\Rcc{p}^{(n)}_{n-k}\diff v_k-\sum_{k=1}^n\Rcc{p}^{(n)}_{n-k}\diff v_k^2
	=\sum_{k=1}^n\Rcc{p}^{(n)}_{n-k}\brab{\diff v_k}\brab{2v_n-v_k-v_{k-1}}\nonumber
\end{align}
and, by using the identity
$2v_n-v_k-v_{k-1}=-(v_k-v_{k-1})+2\sum_{j=k}^n(v_j-v_{j-1})$,
\begin{align}\label{lemp: difference term Jn}
	J_n&=-\sum_{k=1}^n\Rcc{p}^{(n)}_{n-k}\brat{\diff v_k}^2
	+2\sum_{k=1}^n\Rcc{p}^{(n)}_{n-k}\sum_{j=k}^n\brat{\diff v_k}\brat{\diff v_j}\nonumber\\
	&=-\sum_{k=1}^n\Rcc{p}^{(n)}_{n-k}\brat{\diff v_k}^2
	+2\sum_{j=1}^n\brat{\diff v_j}\sum_{k=1}^{j}\Rcc{p}^{(n)}_{n-k}\brat{\diff v_k}.
\end{align}
Define $u_j:=\sum_{k=1}^{j}\Rcc{p}^{(n)}_{n-k}\brat{\diff v_k}$ with $u_0=0$ such that 
$$u_j-u_{j-1}=\Rcc{p}^{(n)}_{n-j}\brat{\diff v_j}\quad\text{and}\quad \diff v_j=\frac{u_j-u_{j-1}}{\Rcc{p}^{(n)}_{n-j}}.$$ 
It follows from \eqref{lemp: difference term Jn} that
\begin{align}\label{lemp: difference term Jn2}
	J_n&=-\sum_{k=1}^n\frac{\brat{u_k-u_{k-1}}^2}{\Rcc{p}^{(n)}_{n-k}}
	+2\sum_{j=1}^n\frac{u_j(u_j-u_{j-1})}{\Rcc{p}^{(n)}_{n-j}}
	=\sum_{j=1}^n\frac{u_j^2-u_{j-1}^2}{\Rcc{p}^{(n)}_{n-j}}\nonumber\\
	&=\frac{u_n^2}{\Rcc{p}^{(n)}_{0}}
	+\sum_{j=1}^n\braB{\frac{1}{\Rcc{p}^{(n)}_{n-j}}-\frac{1}{\Rcc{p}^{(n)}_{n-j-1}}}u_j^2
	=a_{0}^{(n)}u_n^2+\sum_{j=1}^n\braB{\frac{1}{\Rcc{p}^{(n)}_{n-j}}-\frac{1}{\Rcc{p}^{(n)}_{n-j-1}}}u_j^2,
\end{align}
where we use the fact $\Rcc{p}_{0}^{(n)}=\theta_{0}^{(n)}=1/a_{0}^{(n)}$ according to 
\eqref{def: general DOC kernels} and \eqref{def: general RCC kernels}.
By using the equality \eqref{def: wn-vk}, 
we have $u_n=\sum_{k=1}^{n}\Rcc{p}^{(n)}_{n-k}\brat{\diff v_k}=w_n$. 
It follows from \eqref{lemp: difference term Jn0} and  \eqref{lemp: difference term Jn2} that
	\begin{align*}
		2w_n\sum_{j=1}^na_{n-j}^{(n)}w_j=&\,\sum_{k=1}^{n}\Rcc{p}_{n-k}^{(n)}\diff v_k^2
		+a_{0}^{(n)}w_n^2
		+\sum_{j=1}^{n-1}\braB{\frac{1}{\Rcc{p}_{n-j}^{(n)}}-\frac{1}{\Rcc{p}_{n-j-1}^{(n)}}}
		\braB{\sum_{k=1}^{j}\Rcc{p}_{n-k}^{(n)}\diff v_k}^2.
	\end{align*}
Recalling the definition of the auxiliary kernels $a_{n-j}^{(n)}$, one has
	\begin{align*}
		2w_n\sum_{j=1}^{n}\chi_{n-j}^{(n)}w_j
		&=2w_n\sum_{j=1}^{n}a_{n-j}^{(n)}w_j+2\brab{\chi_{0}^{(n)}-a_{0}^{(n)}}w_n^2\nonumber\\
		&=\sum_{k=1}^{n}\Rcc{p}_{n-k}^{(n)}\diff v_k^2
		+\sigma_{\min}\chi_{0}^{(n)}w_n^2
		+\sum_{j=1}^{n-1}\braB{\frac{1}{\Rcc{p}_{n-j}^{(n)}}-\frac{1}{\Rcc{p}_{n-j-1}^{(n)}}}
		\braB{\sum_{k=1}^{j}\Rcc{p}_{n-k}^{(n)}\diff v_k}^2.
	\end{align*}
	Thus the equality \eqref{def: vk-vk2} completes the proof.
\end{proof}

The parameter $\sigma_{\min}$ in Theorem \ref{thm: DGS-general kernels} is set 
	to estimate the minimum eigenvalue of the associated quadratic form with the discrete kernels $\chi_{n-j}^{(n)}$. 
	If the first kernel $\chi_{0}^{(n)}$ is properly large, one can choose a $\sigma_{\min}\in(0,2)$ 
	to meet our assumptions and find a lower bound $\sigma_{\min}\min_k\chi_{0}^{(k)}$ of the minimum eigenvalue,
	see an open problem in Remark \ref{remark: minimum eigenvalue}.
	
\begin{remark}\label{rem: necessity of Lemma 2.1}
To gain better understanding of the parameter $\sigma_{\min}$, we present some further comments to explain
the necessity of Theorem \ref{thm: DGS-general kernels} and its continuous version, 
Lemma \ref{lem: countinous DGS equality}, 
which is an updated version of the following equality in \cite[Lemma 1]{AlsaediAhmadKirane:2015}
\begin{align}\label{RL derivative by AlsaediAhmadKirane:2015}
	2v(t)({}^R\!\partial_t^{\beta}v)(t)
	=\brab{{}^R\!\partial_t^{\beta}v^2}(t)+\omega_{1-\beta}(t)v^2(t)
	-\int_0^t\omega_{-\beta}(t-s)\bra{v(t)-v(s)}^2\zd{s}.	
\end{align}
The main reason for this modification is that the discrete counterpart of \eqref{RL derivative by AlsaediAhmadKirane:2015} may be inadequate to provide a discrete gradient structure for the discrete kernels in \eqref{def: IAF formula kernel} and \eqref{def: L1+ formula kernel}. To see it more clear, we simulate \eqref{RL derivative by AlsaediAhmadKirane:2015} at the discrete time levels by revisiting the difference term $J_n$ in \eqref{lemp: difference term Jn}. 
For the fixed index $n$, let $\Rcc{p}^{(n)}_{n}\equiv0$. One has
$\sum_{k=0}^{j-1}\brab{\Rcc{p}^{(n)}_{n-k-1}-\Rcc{p}^{(n)}_{n-k}}=\Rcc{p}^{(n)}_{n-j},$
and then
\begin{align*}
	2\brat{\diff v_j}\sum_{k=1}^{j}\Rcc{p}^{(n)}_{n-k}\brat{\diff v_k}
	=&2\brat{\diff v_j}\kbra{\Rcc{p}^{(n)}_{n-j}v_{j}-\sum_{k=1}^{j-1}
		\brab{\Rcc{p}^{(n)}_{n-k-1}-\Rcc{p}^{(n)}_{n-k}}v_k-\Rcc{p}^{(n)}_{n-1}v_0}\\
	=&2\sum_{k=0}^{j-1}\brab{\Rcc{p}^{(n)}_{n-k-1}-\Rcc{p}^{(n)}_{n-k}}(v_{j}-v_k)
	\kbra{\brat{v_j-v_{k}}-\brat{v_{j-1}-v_k}}\\
	=&\sum_{k=0}^{j-1}\brab{\Rcc{p}^{(n)}_{n-k-1}-\Rcc{p}^{(n)}_{n-k}}\kbrab{(v_{j}-v_k)^2-(v_{j-1}-v_k)^2}
	+\Rcc{p}^{(n)}_{n-j}(\diff v_j)^2.
\end{align*}
Thus the difference term $J_n$ in \eqref{lemp: difference term Jn} can be handled by
\begin{align*}
	J_n&=-\sum_{k=1}^n\Rcc{p}^{(n)}_{n-k}\brat{\diff v_k}^2
	+2\sum_{j=1}^n\brat{\diff v_j}\sum_{k=1}^{j}\Rcc{p}^{(n)}_{n-k}\brat{\diff v_k}\\
	&=\sum_{j=1}^n\sum_{k=0}^{j-1}\brab{\Rcc{p}^{(n)}_{n-k-1}-\Rcc{p}^{(n)}_{n-k}}
	\kbrab{(v_{j}-v_{k})^2-(v_{j-1}-v_{k})^2}\\
	&=\sum_{k=0}^{n-1}\brab{\Rcc{p}^{(n)}_{n-k-1}-\Rcc{p}^{(n)}_{n-k}}
	\sum_{j=k+1}^{n}\kbrab{(v_{j}-v_{k})^2-(v_{j-1}-v_{k})^2}\\
	&=\sum_{k=1}^{n-1}\brab{\Rcc{p}^{(n)}_{n-k-1}-\Rcc{p}^{(n)}_{n-k}}(v_{n}-v_{k})^2
	+\Rcc{p}^{(n)}_{n-1}(v_{n}-v_{0})^2.	
\end{align*}
By using the definition \eqref{lemp: difference term Jn0} of $J_n$ 
and the definition of $a_{n-j}^{(n)}$, it is easy to obtain that
\begin{align*}
	2w_n\sum_{j=1}^na_{n-j}^{(n)}w_j=&\,\sum_{k=1}^{n}\Rcc{p}_{n-k}^{(n)}\diff v_k^2
	+\sum_{k=1}^{n-1}\brab{\Rcc{p}^{(n)}_{n-k-1}-\Rcc{p}^{(n)}_{n-k}}(v_{n}-v_{k})^2
	+\Rcc{p}^{(n)}_{n-1}(v_{n}-v_{0})^2
\end{align*}
and thus
\begin{align*}
	2w_n\sum_{j=1}^{n}\chi_{n-j}^{(n)}w_j
	=&\,2(\sigma_{\min}-1)\chi_{0}^{(n)}w_n^2
	+\sum_{k=1}^{n}p_{n-k}^{(n)}v_k^2-\sum_{k=1}^{n-1}p_{n-k-1}^{(n-1)}v_k^2\\
	&\,+\sum_{k=1}^{n-1}\brab{\Rcc{p}^{(n)}_{n-k-1}-\Rcc{p}^{(n)}_{n-k}}(v_{n}-v_{k})^2
	+\Rcc{p}^{(n)}_{n-1}(v_{n}-v_{0})^2.
\end{align*}
As seen, a desired DGS requires $\sigma_{\min}\ge1$ 
(that is, the first kernel $\chi_{0}^{(n)}$ has to be properly large);
nonetheless, this requirement can not be fulfilled for
the discrete kernels in \eqref{def: IAF formula kernel} and \eqref{def: L1+ formula kernel}.
As we will see in subsequent discussions, a small parameter $\sigma_{\min}\in[0,1-\beta]$ is always necessary 
since the first kernel $a_{0}^{(\beta,n)}$ is not always dominant.
\end{remark}

\subsection{DGS of integral averaged formula}

To build up the desired DGS \eqref{eq: DGS equality},
it remains to verify that the kernels $a_{n-k}^{(\beta,n)}$ in \eqref{def: IAF formula kernel}
meet the assumptions of Theorem \ref{thm: DGS-general kernels}.
With the definition \eqref{def: IAF formula kernel}, the integral mean-value theorem yields
	\[
	a_{0}^{(\beta,n)}=\frac{\tau_{n}^{\beta-1}}{\Gamma(2+\beta)}\quad\text{and}\quad
	a_{1}^{(\beta,n)}>a_{2}^{(\beta,n)}>\cdots>a_{n-1}^{(\beta,n)}>0\quad \text{for $n\ge2$}.
	\]
A direct calculation gives
\begin{align*}
	a_{0}^{(\beta,n)}-a_{1}^{(\beta,n)}
	=\frac{r_n}{\Gamma(2+\beta)\tau_{n}^{1-\beta}}
	\kbra{1+1/r_n+1/r_n^{1+\beta}-(1+1/r_n)^{1+\beta}}\quad\text{for $n\ge2$.}
\end{align*}
It is easily seen that $a_{0}^{(\beta,n)}>a_{1}^{(\beta,n)}$ as $\beta\rightarrow{0}$,
while $a_{0}^{(\beta,n)}<a_{1}^{(\beta,n)}$ as $\beta\rightarrow{1}$.
The value of $a_{0}^{(\beta,n)}-a_{1}^{(\beta,n)}$ changes the sign when
the fractional index $\beta$ varies over $(0,1)$.
The discrete kernels $a_{j}^{(\beta,n)}$ are not uniformly monotonous with respect to the subscript $j$
so that the recent theory \cite{JiZhuLiao:2022,LiaoMcLeanZhang:2019,LiaoTangZhou:2020,LiaoZhuWang:2022NMTMA,StynesORiordanGracia:2017}
for the nonuniform L1 and L2-1$_{\sigma}$ formulas can not be directly
applied to the numerical analysis of the integral averaged formulas \eqref{def: IAF formula}
and \eqref{def: L1+ formula}.

On the other hand, it is easy to check that $(1+x)^\beta<1+x^\beta$, which leads to
\begin{align*}
	(1+x)^{1+\beta}<1+\brat{1+\beta}x+x^{1+\beta}\quad\text{for $x>0$.}
\end{align*}
We see that
\begin{align*}
	(1+\beta)a_{0}^{(\beta,n)}-a_{1}^{(\beta,n)}
	=\frac{r_n}{\Gamma(2+\beta)\tau_{n}^{1-\beta}}
	\kbra{1+(1+\beta)/r_n+1/r_n^{1+\beta}-(1+1/r_n)^{1+\beta}}>0
\end{align*}
for $n\ge2$. In summary, one has the following result.
\begin{lemma}\label{lem: L1+ Coeff-decreasing}
	The discrete kernels $a_{n-k}^{(\beta,n)}$ in \eqref{def: IAF formula kernel} fulfill
	\[
	(1+\beta)a_{0}^{(\beta,n)}>a_{1}^{(\beta,n)}>a_{2}^{(\beta,n)}>\cdots>a_{n-1}^{(\beta,n)}>0
	\quad \text{for $n\ge2$}.
	\]
\end{lemma}

To process the analysis, we define a function 
\begin{align}\label{def: app-varphi}
	\rho(z):=(z+1)^{1+\beta}-z^{1+\beta}-1\quad\text{for $z\ge0$}.
\end{align}

\begin{lemma}\label{lem: L1+ geometric convexity}
Let the adjacent step-ratios satisfy the following condition
	\begin{align}\label{ieq: condition}
		r_{k+1}\geq r_*(r_{k}):=\kbra{\frac{(2^\beta-1)\rho(r_{k})}{\rho(2r_{k})-\rho(r_{k})}}^{\frac1{1-\beta}}\quad\text{for $k\ge2$},
	\end{align}
where the function $\rho$ is defined in \eqref{def: app-varphi} and $r_*(z)<1$ for any $z>0$.	
	Then the discrete convolution kernels $a_{n-k}^{(\beta,n)}$ in \eqref{def: IAF formula kernel} fulfill
	\begin{align}\label{ieq: L1+ geometric convexity}
	\frac{a^{(\beta,n)}_{1}}{2a^{(\beta,n-1)}_{0}}<
    \frac{a^{(\beta,n)}_{2}}{a^{(\beta,n-1)}_{1}}
	<\cdots<
    \frac{a^{(\beta,n)}_{n-1}}{a^{(\beta,n-1)}_{n-2}}<1\quad\text{for $n\ge2$.} 	
\end{align}
    \end{lemma}

\begin{proof}	
	We define a class of auxiliary function
	\begin{align*}
		a_{n,k}(\mu)&:=\frac{1}{\tau_n\tau_{k}}\int_{t_{n-1}}^{t_n}\int_{t_{k-1}}^{t_{k-1}+\mu\tau_k}
		\omega_{\beta}(t-s)\zd{s}\zd{t}\quad\text{for $1\leq{k}\leq{n-1}$},
	\end{align*}
	such that $a^{(\beta,n)}_{n-k}=a_{n,k}(1)$ for $1\leq{k}\leq{n-1}$, $a_{n,k}(0)=0$ and
	\begin{align*}
		a_{n,k}'(\mu)=\frac{1}{\tau_n}\int_{t_{n-1}}^{t_{n}}\omega_{\beta}(t-t_{k-1}-\mu\tau_k)\zd{t}
		\quad\text{for $1\leq{k}\leq{n-1}$}.
	\end{align*}
	Consider the following sequence
	\begin{align}\label{lemproof: CG auxiliary psi}
		\psi_{1}^{(n)}:=\frac{a^{(\beta,n)}_{1}}{2a^{(\beta,n-1)}_{0}}
		\quad\text{and}\quad\psi_{n-k}^{(n)}:=\frac{a^{(\beta,n)}_{n-k}}{a^{(\beta,n-1)}_{n-k-1}}
		=\frac{a_{n,k}(1)}{a_{n-1,k}(1)}
		\quad \text{for $1\le k\le n-2$.}
	\end{align}
	Thanks to the Cauchy differential
	mean-value theorem, there exists $\mu_{nk}\in(0,1)$ such that
	\begin{align}\label{lemproof: CG psi equality}
		\psi_{n-k}^{(n)}&=\frac{a_{n,k}(1)}{a_{n-1,k}(1)}=
		\frac{a_{n,k}(1)-a_{n,k}(0)}{a_{n-1,k}(1)-a_{n-1,k}(0)}
		=\frac{a_{n,k}'(\mu_{nk})}{a_{n-1,k}'(\mu_{nk})}\nonumber\\
		&=\frac{\int_{t_{n-1}}^{t_{n}}\omega_{\beta}(t-t_{k-1}-\mu_{nk}\tau_k)\zd{t}}
		{r_n\int_{t_{n-2}}^{t_{n-1}}\omega_{\beta}(t-t_{k-1}-\mu_{nk}\tau_k)\zd{t}}
		\quad \text{for $1\leq k\leq n-2$}.
	\end{align}
	It is not difficult to verify 
	that the following function  
	\begin{align}\label{lemproof: CG auxiliary gn}
		g_n(z):=&\,\frac{\int_{t_{n-1}}^{t_{n}}\omega_{\beta}(t-z)\zd{t}}
		{r_n\int_{t_{n-2}}^{t_{n-1}}\omega_{\beta}(t-z)\zd{t}}\quad\text{for fixed $n\ge3$}
	\end{align}
is decreasing with respect to $z\in(0, t_{n-2})$ by checking the sign of the first derivative $g_n'(z)$.
	Thus we derive from \eqref{lemproof: CG psi equality} that 
	\begin{align*}
		g_n(t_k)<\psi_{n-k}^{(n)}<g_n(t_{k-1})
		\quad \text{for $1\leq k\leq n-2$}.
	\end{align*}
	By taking $k=n-2,n-3,\cdots,1$ in this inequality, we arrive at
	\begin{align}\label{lemproof: CG psi inequality}
		g_n(t_{n-2})<\psi_{2}^{(n)}<\psi_{3}^{(n)}<\cdots<\psi_{n-1}^{(n)}<g_n(t_0)<1
		\quad \text{for $n\ge3$},
	\end{align}
	where we use the simple fact $g_n(t_0)<1$ due to the integral mean-value theorem.
	
	It remains to check that $\psi_{1}^{(n)}<\psi_{2}^{(n)}$ for $n\ge3$.
	Applying \eqref{def: IAF formula kernel} and \eqref{lemproof: CG auxiliary psi}, we have
	\begin{align*}
		\psi_{1}^{(n)}=&\,\frac{(r_{n}+1)^{1+\beta}-r_{n}^{1+\beta}-1}{2r_n}=\frac{\rho(r_{n})}{2r_n},\\
		\psi_2^{(n)}
		=&\,\frac{(r_n r_{n-1}+r_{n-1}+1)^{1+\beta}-(r_n r_{n-1}+r_{n-1})^{1+\beta}
			-(r_{n-1}+1)^{1+\beta}+r_{n-1}^{1+\beta}}{r_n\kbrab{(r_{n-1}+1)^{1+\beta}-r_{n-1}^{1+\beta}-1}}\\
		=&\,\frac{\rho(r_{n}r_{n-1}+r_{n-1})-\rho(r_{n-1})}{r_n\rho(r_{n-1})}\,.		
	\end{align*}
	Recalling the definition \eqref{lemproof: CG auxiliary gn}, one has
	$g_n(t_{n-2})=\kbrat{(r_{n}+1)^{\beta}-1}/r_n$. Then Lemma \ref{lem: app-hx} and 
	the first inequality of \eqref{lemproof: CG psi inequality} yield
	\begin{align}\label{lemproof: CG psi1-decrese1}
		\psi_{1}^{(n)}=\frac{\rho(r_{n})}{2r_n}\le g_n(t_{n-2})<\psi_{2}^{(n)}
		\quad \text{for $n\ge3$\; if $r_n\ge1$.}
	\end{align}
	Specially, taking $r_n=1$ in \eqref{lemproof: CG psi1-decrese1} gives
	\begin{align*}
		\frac{\rho(2r_{n-1})-\rho(r_{n-1})}{\rho(r_{n-1})}> \frac{\rho(1)}{2}=2^\beta-1
		\quad \text{for $n\ge3$.}
	\end{align*}
Then the definition \eqref{ieq: condition} of $r_*$ with the arbitrariness of  $r_{n-1}$ implies that 
\begin{align}\label{ieq: r* condition}
	r_{*}(z)<1
	\quad \text{for any $z>0$.}
\end{align}
Furthermore, with the adjacent step-ratio constraint \eqref{ieq: condition}, Lemma \ref{lem: app-Hxy} shows that
	\begin{align}\label{lemproof: CG psi1-decrese2}
		\psi_{1}^{(n)}<\psi_{2}^{(n)}
		\quad \text{for $n\ge3$\; if $r_{*}(r_{n-1})\le r_n<1$.}
	\end{align}
The desired result \eqref{ieq: L1+ geometric convexity} 
 follows from \eqref{lemproof: CG psi inequality}-\eqref{lemproof: CG psi1-decrese2} immediately.
\end{proof}

\begin{figure}[htb!]
	\centering
		\includegraphics[width=3.4in]{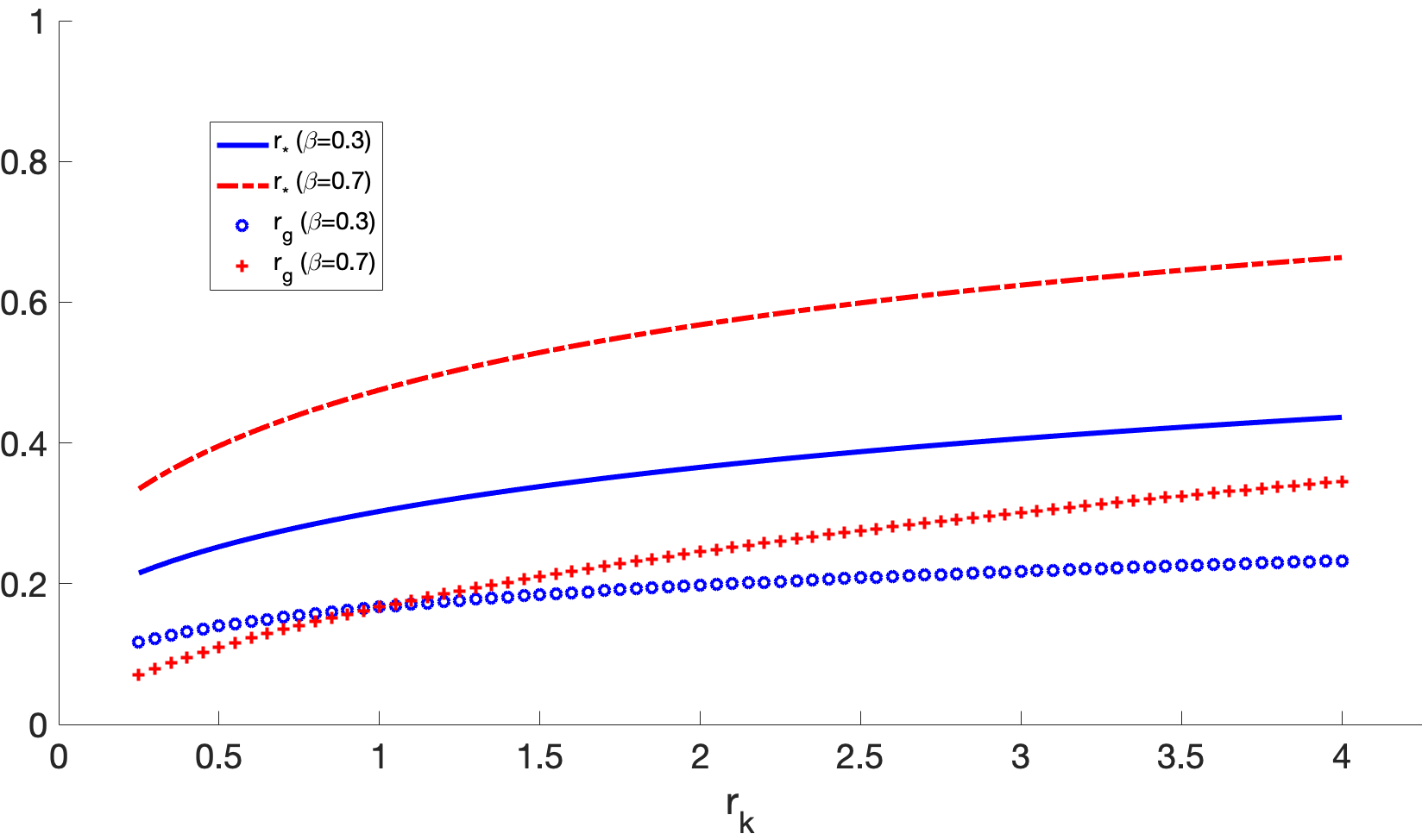}
	\caption{Example curves of $r_*(r_k)$  and $r_g(r_k)$, defined in 
		\eqref{ieq: condition} and \eqref{ieq: guess condition} respectively, 
		for the current step-ratio $r_k\in(1/4,4)$ with two fractional orders $\beta=0.3$ and $0.7$. }
	\label{fig: step-ratio condition}
\end{figure}

	\begin{remark}\label{remark: step-ratio condition}
	It is to remark that the step-ratio constraint \eqref{ieq: condition} 
	for Lemma \ref{lem: L1+ geometric convexity} is imposed rather theoretically than practically.
	For the current time-step size $\tau_{k}$, this step-ratio constraint \eqref{ieq: condition}
	allows the next time-step size $\tau_{k+1}$ to increase freely or decrease within a certain range. 
	As seen in Figure \ref{fig: step-ratio condition},
	the step-ratio constraint \eqref{ieq: condition} is practically mild in numerical simulations although 
	it is only a sufficient condition.  Our numerical computations show that
	a weak restriction 
	\begin{align}\label{ieq: guess condition}
		r_{k+1}\geq r_{g}(r_k):=\brab{1+5 r_{k}^{-\beta}}^{-1}
		\quad\text{for $k\ge2$}
	\end{align}
	is also sufficient to ensure the desired inequality \eqref{lemproof: CG psi1-decrese2}
	or Lemma \ref{lem: L1+ geometric convexity}; 
	however, we are not able to present a rigorous proof under 
	this updated restriction \eqref{ieq: guess condition}.
\end{remark}

Lemmas \ref{lem: L1+ Coeff-decreasing} and \ref{lem: L1+ geometric convexity}
say that the discrete convolution kernels $a_{n-k}^{(\beta,n)}$ in \eqref{def: IAF formula kernel}
meet the assumptions of Theorem \ref{thm: DGS-general kernels} with $\sigma_{\min}=0$.
Thus we have the following theorem.
\begin{theorem}\label{thm: DGS-IAF kernels}
	If the step-ratios satisfy the condition \eqref{ieq: condition},
	the DGS equality \eqref{eq: DGS equality} holds and the discrete kernels $a_{n-k}^{(\beta,n)}$
	in \eqref{def: IAF formula kernel}
	are positive definite in the sense of \eqref{ieq: positive definite}.	
\end{theorem}

\begin{remark}[The minimum eigenvalue and an open problem]\label{remark: minimum eigenvalue}	
	The imposed step-ratio constraint \eqref{ieq: condition} and the main technical difficulty
	 (see Lemmas \ref{lem: app-hx}-\ref{lem: app-Hxy}) in Lemma \ref{lem: L1+ geometric convexity} 
	are raised for establishing the first inequality of \eqref{ieq: L1+ geometric convexity},
	while the other inequalities of \eqref{ieq: L1+ geometric convexity} are valid 
	without any step-ratio conditions on arbitrary time meshes. 
	By numerical computations (no theoretical proof), one may find an implicit function
	 $r_{E}=r_{E}(r_k)$ such that
		\[
	\frac{a^{(\beta,n)}_{1}}{(1+\beta)a^{(\beta,n-1)}_{0}}<
	\frac{a^{(\beta,n)}_{2}}{a^{(\beta,n-1)}_{1}}
	\quad\text{for $n\ge2$\quad if $r_{k+1}\ge r_{E}(r_k)$ for $k\ge2$.} \]  
	In this situation, the discrete kernels $a_{n-k}^{(\beta,n)}$ 
	in \eqref{def: IAF formula kernel}
	will meet the assumptions of Theorem \ref{thm: DGS-general kernels} 
	with $\sigma_{\min}=1-\beta$. Then the inequality \eqref{ieq: positive definite} can be updated by
	\begin{align*}
		2\sum_{k=1}^nw_k\sum_{j=1}^ka^{(\beta,k)}_{k-j}w_j
		\ge (1-\beta)\sum_{k=1}^{n}a_{n-k}^{(\beta,n)}w_k^2
		+\sum_{k=1}^{n}\mathsf{p}_{n-k}^{(\beta,n)}v_k^2\quad\text{for nonzero $\{w_{k}\}_{k=1}^n$.}
	\end{align*}
	\begin{table}[!ht]
	\begin{center}
		\tabcolsep 0pt {Table 1 \quad The minimum eigenvalue $\lambda_{\min}$ 
			and the bound $\sigma_\beta$
			on random mesh.} \vspace*{0.5pt}
		\def\temptablewidth{0.9\textwidth}
		{\rule{\temptablewidth}{1pt}}
		\begin{tabular*}{\temptablewidth}{@{\extracolsep{\fill}}ccccccc}
			$n$ &\multicolumn{2}{c}{$\beta=0.8$}&\multicolumn{2}{c}{$\beta=0.5$}
			&\multicolumn{2}{c}{$\beta=0.1$} \\
			\cline{2-3}    \cline{4-5} \cline{6-7}
			& $\lambda_{\min}$  & $\sigma_\beta$ &$\lambda_{\min}$  
			& $\sigma_\beta$ &$\lambda_{\min}$ &$\sigma_\beta$\\
			\hline
			$100$   &0.4999  & 0.2599  &5.7501   &2.7519 &59.2896   &27.8949 \\
			$200$   &0.5842  &0.3038   &7.6961   &3.6778 &118.5504 &56.2933 \\
			$400$   &0.6676  &0.3445   &11.3192 &5.4621 &202.7272 &96.3790 \\			
		\end{tabular*}
		{\rule{\temptablewidth}{0.9pt}}
	\end{center}\label{table: random minimum eigenvalue}
\end{table}
\begin{table}[!ht]
	\begin{center}
		\tabcolsep 0pt {Table 2 \quad The minimum eigenvalue $\lambda_{\min}$ 
			and the bound $\sigma_{1/2}$
			on graded mesh.} \vspace*{0.5pt}
		\def\temptablewidth{0.9\textwidth}
		{\rule{\temptablewidth}{1pt}}
		\begin{tabular*}{\temptablewidth}{@{\extracolsep{\fill}}ccccccc}
			$n$ &\multicolumn{2}{c}{$\gamma=1.5$}&\multicolumn{2}{c}{$\gamma=2.0$}
			&\multicolumn{2}{c}{$\gamma=3.0$} \\
			\cline{2-3}    \cline{4-5} \cline{6-7}
			& $\lambda_{\min}$  & $\sigma_{1/2}$ &$\lambda_{\min}$  
			& $\sigma_{1/2}$ &$\lambda_{\min}$ &$\sigma_{1/2}$\\
			\hline
			$100$   &5.9665   & 3.0749 &5.2342   &2.6663  &4.3604 &2.1825 \\
			$200$   &8.3667   &4.3459  &7.3055   &3.7660  &6.0420 &3.0788 \\
			$400$   &11.7695 &6.1440  &10.2459 &5.3226  &8.4345 &4.3486 \\			
		\end{tabular*}
		{\rule{\temptablewidth}{0.9pt}}
	\end{center}\label{table: graded minimum eigenvalue}
\end{table}
	It suggests that the quadratic form with the discrete kernels $a_{n-k}^{(\beta,n)}$ 
	has a minimum eigenvalue 
	larger than 
	$\sigma_\beta:=(1-\beta)\min_k a_{0}^{(\beta,k)}$. 
	For three different fractional orders $\beta=0.8,0.5$ and $0.1$, 
	 Table 1 lists the minimum eigenvalue and the lower bound $\sigma_\beta$  	
	of the associated quadratic form on random time meshes $(t_n=1)$.
	With  the fractional order $\beta=1/2$, the data in Table 2 is computed 
	on the graded mesh $t_k=(k/n)^{\gamma}$ with three different grading parameters $\gamma=1.5,2.0$ and $3.0$.  
	These numerical results support our theoretical prediction although we are not able to verify it theoretically.  	
\end{remark}

\section{Application to time-fractional Allen-Cahn model}
\setcounter{equation}{0}

In this section, we consider
the time-fractional Allen-Cahn model \cite{LiaoTangZhou:2020,LiaoTangZhou:2021, LiaoZhuWang:2022NMTMA,QuanTangYang-2020csiam,
	QuanTangYang-2020,QuanTangWangYang:2022,TangYuZhou:2019}
\begin{align}\label{cont: TFAC model}
\partial_t^\alpha \Phi=-\kappa\mu\quad\text{with
the chemical potential $\mu:=\frac{\delta E}{\delta \Phi}=f(\Phi)-\epsilon^2\Delta\Phi$},
\end{align}
where $\kappa$ is the mobility coefficient and
$E$ is the Ginzburg-Landau energy functional
\begin{align}\label{cont: classical free energy}
E[\Phi]:=\int_{\Omega}\braB{\frac{\epsilon^2}{2}\abs{\nabla\Phi}^2
+F(\Phi)}\zd\mathbf{x}\quad\text{with\; $F(\Phi):=\frac{1}{4}\bra{\Phi^2-1}^2$.}
\end{align}
Here,  the real valued function $\Phi$ represents the concentration
difference in a binary system,
$\epsilon>0$ is an interfacial width parameter and the nonlinear term $f(\Phi)=F'(\Phi)$.
For simplicity, consider the spatial domain $\Omega\subseteq \mathbb{R}^2$
with the periodic boundary conditions.

\subsection{Continuous energy dissipation law}

Let $\myinner{\cdot,\cdot}$ and $\mynorm{\cdot}$ denote the
$L^2(\Omega)$ inner product and the associated norm, respectively.
Also, we use the standard norms of the Sobolev space $H^{m}\bra{\Omega}$
and the $L^p\bra{\Omega}$ space.
By the model \eqref{cont: TFAC model},
the time derivative of free energy
\begin{align}\label{cont: TFAC derivative-energy}
	\frac{\zd }{\zd t}E[\Phi]=\myinnerb{\tfrac{\delta E}{\delta \Phi},\partial_t\Phi}
	=-\frac1{\kappa}\myinnerb{\partial_t^\alpha \Phi,\partial_t\Phi}.
\end{align}
Taking $\beta:=1-\alpha$ and $w:=\partial_t\Phi$ in Lemma \ref{lem: countinous DGS equality}
with
$$v=\mathcal{I}_t^{1-\alpha}w=\partial_t^\alpha \Phi=-\kappa\mu,$$
one gets the following equality
\begin{align*}
\frac1{\kappa}\myinnerb{\partial_t^\alpha \Phi,\partial_t\Phi}
	=&\,\frac{\kappa}2\frac{\zd }{\zd t}\brab{\mathcal{I}_{t}^{\alpha}\mynormb{\mu}^2}
	+\frac{1}{2\kappa}
	\int_0^t\frac{\partial}{\partial\xi}\bra{\frac{1}{\omega_{\alpha}(t-\xi)}}
	\mynormB{\int_0^{\xi}\omega_{\alpha}(t-s)v'(s)\zd{s}}^2\zd \xi.
\end{align*}
Thus by using \eqref{cont: TFAC derivative-energy} we get a variational energy dissipation law
\begin{align}\label{cont: TFAC energy law}
	\frac{\zd E_{\alpha}}{\zd t}+\frac{1}{2\kappa}
	\int_0^t\frac{\partial}{\partial\xi}\bra{\frac{1}{\omega_{\alpha}(t-\xi)}}
	\mynormB{\int_0^{\xi}\omega_{\alpha}(t-s)v'(s)\zd{s}}^2\zd \xi=0,
\end{align}
where the nonlocal (variational) energy
\begin{align}\label{cont: TFAC modified energy}
	E_{\alpha}[\Phi]:=E[\Phi]
	+\frac{\kappa}{2}\mathcal{I}_{t}^{\alpha}\mynormb{\mu}^2=E[\Phi]
	+\frac{\kappa}{2}\mathcal{I}_{t}^{\alpha}\mynormb{\tfrac{\delta E}{\delta \Phi}}^2\quad\text{for $t>0$.}
\end{align}
This new energy dissipation law updates the previous energy laws in \cite{JiZhuLiao:2022,LiaoTangZhou:2021,LiaoZhuWang:2022NMTMA}
in the sense that it can exactly recover the energy dissipation law of the classical Allen-Chan equation.
As the fractional order $\alpha\rightarrow1$,
one can check that
$\int_0^{\xi}\omega_{\alpha}(t-s)v'(s)\zd{s}\rightarrow v(\xi)=-\kappa \mu(\xi)$
 and the variational energy dissipation law \eqref{cont: TFAC energy law} degrades into
\begin{align*}
	\frac{\zd }{\zd t}\brab{E[\Phi]
		+\frac{\kappa}{2}\mathcal{I}_{t}^{1}\mynormb{\mu}^2}+\frac{\kappa}{2}\mynormb{\mu}^2
	=\frac{\zd E}{\zd t}+\kappa\mynormb{\mu}^2=0,
\end{align*}
which is just the energy dissipation law of the classical Allen-Cahn model.
In this sense, we say that the energy law \eqref{cont: TFAC energy law}
is asymptotically compatible in the  limit $\alpha\rightarrow1$.

\subsection{Crank-Nicolson scheme}

We will only consider the time-discrete methods, here and hereafter, with the numerical solution
$\phi^n\approx \Phi(t_n)$.
Our numerical scheme and the analysis can be extended in a straightforward way to the
fully discrete schemes with some appropriate spatial discretization
preserving the discrete integration-by-parts formulas.
By applying the L1$^+$ formula \eqref{def: L1+ formula}, we have the  following Crank-Nicolson scheme
\begin{align}\label{scheme: CN TFAC model}
	\bra{\partial_{\tau}^\alpha \phi}^{n-\frac12}
	=-\kappa\mu^{n-\frac12}\quad\text{with}\quad
	\mu^{n-\frac12}=f\brat{\phi}^{n-\frac12}-\epsilon^2\Delta \phi^{n-\frac12}
	\quad\text{for $n\ge 1$.}
\end{align}
Here, $f\brat{\phi}^{n-\frac12}$ is the standard second-order
approximation defined by
\begin{align}\label{def: nonlinear approximation}
f\brat{\phi}^{n-\frac12}:=\frac12\kbrab{\brat{\phi^n}^2+\brat{\phi^{n-1}}^2}\phi^{n-\frac12}-\phi^{n-\frac12}
\end{align}
such that
\begin{align}\label{eq: property nonlinear approximation}
	\myinnerb{f\brat{\phi}^{n-\frac12},\diff \phi^n}
=\myinnerb{F(\phi^n),1}-\myinnerb{F(\phi^{n-1}),1}.
\end{align}

Recalling the L1$^+$ kernels $a^{(1-\alpha,n)}_{n-k}$ defined in \eqref{def: L1+ formula}
with 
$L^{n-1}:=\sum_{k=1}^{n-1}a^{(1-\alpha,n)}_{n-k}\diff \phi^{k}$,
we consider the following discrete functional $G[z]$,
\begin{align*}
G[z]:=&\,\frac{1}{2\kappa}a^{(1-\alpha,n)}_{0}\mynormb{z-\phi^{n-1}}^2
+\frac{1}{\kappa}\myinnerb{L^{n-1},z-\phi^{n-1}}
+\frac{\epsilon^2}{4}\mynormb{\nabla \brab{z+\phi^{n-1}}}^2\\
&\,+\frac{1}{16}\mynormb{z}_{L^4}^4
+\frac{1}{12}\myinnerb{\phi^{n-1},z^3}
+\frac{1}{8}\myinnerb{\brat{\phi^{n-1}}^2,z^2}
+\frac{1}{4}\myinnerb{\brat{\phi^{n-1}}^3,z}
-\frac{1}4\mynormb{z+\phi^{n-1}}^2.
\end{align*}
The solution $\phi^n$ of nonlinear equation \eqref{scheme: CN TFAC model} at the time level $t_n$
 is equivalent to the minimum of $G[z]$ if and only if it is strictly convex and coercive.
With the requirement $a^{(1-\alpha,n)}_{0}\ge \kappa/2$ for the convexity,
one can follow the proof of \cite[Theorem 2.1]{LiaoZhuWang:2022NMTMA} to prove the following result.

\begin{lemma}
If $\tau_n\le \sqrt[\alpha]{\frac{2}{\kappa\Gamma(3-\alpha)}}$,
the Crank-Nicolson scheme \eqref{scheme: CN TFAC model} is uniquely solvable.
\end{lemma}

\subsection{Discrete energy dissipation law}

By virtues of Theorem \ref{thm: DGS-IAF kernels},
we establish  a discrete energy dissipation law for
the Crank-Nicolson scheme \eqref{scheme: CN TFAC model}.
With the original energy $E\kbra{\phi^{n}}$
defined via \eqref{cont: classical free energy}, we define the following discrete variational energy
\begin{align}\label{def: variational energy TFAC}
E_\alpha\kbrat{\phi^{n}}:=E\kbrat{\phi^{n}}
+\frac{1}{2\kappa}\sum_{j=1}^n\mathsf{p}_{n-j}^{(1-\alpha,n)}\mynormb{v^j}^2
\quad\text{with $v^j:=\sum_{\ell=1}^{j}\mathsf{a}_{j-\ell}^{(1-\alpha,j)}\diff\phi^{\ell}$.}
\end{align}

\begin{theorem}\label{thm: CN energy law TFAC}
Under the step-ratio constraint \eqref{ieq: condition}, the 
variable-step Crank-Nicolson scheme \eqref{scheme: CN TFAC model}
is energy stable in the sense that it
preserves a discrete energy dissipation law
\begin{align*}
\partial_{\tau}E_\alpha\kbrat{\phi^n}
+\frac{1}{2\kappa\tau_n}\sum_{j=1}^{n-1}\braB{\frac{1}{\Rcc{\mathsf{p}}_{n-j}^{(1-\alpha,n)}}
	-\frac{1}{\Rcc{\mathsf{p}}_{n-j-1}^{(1-\alpha,n)}}}
\mynormB{\sum_{k=1}^{j}\Rcc{\mathsf{p}}_{n-k}^{(1-\alpha,n)}\diff v^k}^2=0.
\end{align*}
\end{theorem}
\begin{proof}
Taking the inner product of  \eqref{scheme: CN TFAC model} by  $\diff \phi^{n}/\kappa$, one gets
\begin{align*}
\frac{1}{\kappa}\myinnerb{\bra{\partial_{\tau}^\alpha\phi}^{n-\frac12},
\diff \phi^n}
+\myinnerb{f\brat{\phi}^{n-\frac12},\diff \phi^n}
+\epsilon^2\myinnerb{\nabla\phi^{n-\frac12},\diff \nabla\phi^n}=0.
\end{align*}
By using the equality \eqref{eq: property nonlinear approximation}, we get
\begin{align}\label{thmProof:L1+ energy law inner}
	\frac{1}{\kappa}\myinnerb{\bra{\partial_{\tau}^\alpha\phi}^{n-\frac12},
		\diff \phi^n}
	+E\kbrat{\phi^n}-E\kbrat{\phi^{n-1}}=0.
\end{align}
Taking $\beta:=1-\alpha$ and $w_n:=\diff\phi^n$ in \eqref{eq: DGS equality}, we have
\begin{align*}
	\frac{1}{\kappa}\myinnerb{\bra{\partial_{\tau}^\alpha\phi}^{n-\frac12},	\diff \phi^n}=&\,
	\frac{1}{2\kappa}\sum_{j=1}^{n}\mathsf{p}_{n-j}^{(1-\alpha,n)}\mynormb{v^j}^2
	-\frac{1}{2\kappa}\sum_{j=1}^{n-1}\mathsf{p}_{n-1-j}^{(1-\alpha,n-1)}\mynormb{v^j}^2\nonumber\\
	&\,+\frac{1}{2\kappa}\sum_{j=1}^{n-1}\braB{\frac{1}{\Rcc{\mathsf{p}}_{n-j}^{(1-\alpha,n)}}
		-\frac{1}{\Rcc{\mathsf{p}}_{n-j-1}^{(1-\alpha,n)}}}
	\mynormB{\sum_{k=1}^{j}\Rcc{\mathsf{p}}_{n-k}^{(1-\alpha,n)}\diff v^k}^2,
\end{align*}
where the sequence $v^j:=\sum_{\ell=1}^{j}\mathsf{a}_{j-\ell}^{(1-\alpha,j)}\diff\phi^{\ell}$.
Inserting it into \eqref{thmProof:L1+ energy law inner}
yields the discrete energy dissipation law immediately.
This completes the proof.
\end{proof}


As the fractional order $\alpha\rightarrow1$,
the definition \eqref{def: L1+ formula kernel} gives the values
$a_{0}^{(0,n)}=1/\tau_{n}$ and $a_{n-k}^{(0,n)}=0$ for $1\le k\le n-1$.
The L1$^+$ time-stepping scheme \eqref{scheme: CN TFAC model}
degrades into the following Crank-Nicolson scheme
\begin{align}\label{scheme: CN allen-cahn model}
	\partial_{\tau} \phi^{n}
	=-\kappa\mu^{n-\frac12}\quad\text{with}\quad
	\mu^{n-\frac12}=f\brat{\phi}^{n-\frac12}-\epsilon^2\Delta \phi^{n-\frac12}
	\quad\text{for $n\ge 1$.}
\end{align}
It is uniquely solvable if $\tau_n\le 2/{\kappa}$ and preserves the following discrete energy law,
\[
\partial_\tau E\kbra{\phi^n}
+\kappa\mynormb{\mu^{n-\frac12}}^2= 0
\quad \text{for  $n\ge1$.}
\]
The definition \eqref{def: modified kernels} shows that the modified kernels $\mathsf{a}_{0}^{(0,n)}=2/\tau_{n}$ and $\mathsf{a}_{n-k}^{(0,n)}=0$ for $1\le k\le n-1$.
By \eqref{def: general DOC kernels},
the associated DOC kernels $\theta_{0}^{(0,n)}=\tau_{n}/2$ and $\theta_{n-k}^{(0,n)}=0$ for $1\le k\le n-1$.
Then the DCC and RCC kernels
\[
\mathsf{p}_{n-k}^{(0,n)}=\tau_{k}/2\quad \text{and}\quad\Rcc{\mathsf{p}}_{n-k}^{(0,n)}=\tau_{n}/2
\quad \text{for $1\le k\le n$}.
\]
Then with $w_k:=\diff\phi^k$ and $v_j:=2\diff \phi^{j}/\tau_j$, the DGS equality \eqref{eq: DGS equality} degrades into
\begin{align*}
	\diff\phi^k\sum_{j=1}^ka^{(0,k)}_{k-j}\diff\phi^j=\frac1{4}\sum_{j=1}^k\tau_jv_j^2
	-\frac1{4}\sum_{j=1}^{k-1}\tau_jv_j^2=\frac{1}{\tau_k}\brat{\diff \phi^{k}}^2\quad\text{for $k\ge1$.}
\end{align*}
That is, the DGS equality \eqref{eq: DGS equality} is asymptotically compatible in the limit $\alpha\rightarrow 1.$
Obviously, the above discrete variational energy \eqref{def: variational energy TFAC} degrades into
\begin{align*}
	E_\alpha\kbrat{\phi^{n}} \quad\longrightarrow\quad E\kbrat{\phi^{n}}
	+\frac{1}{\kappa}\sum_{j=1}^n\tau_j\mynormb{\partial_{\tau}\phi^j}^2\quad\text{as $\alpha\rightarrow 1$;}
\end{align*}
and the discrete energy dissipation law in Theorem \ref{thm: CN energy law TFAC} degrades into
\[
\partial_\tau E\kbra{\phi^n}
+\frac1{\kappa}\mynormb{\partial_{\tau}\phi^n}^2= 0
\quad \text{for  $n\ge1$,}
\]
which is just the discrete energy law of \eqref{scheme: CN allen-cahn model} since $\partial_{\tau}\phi^n=-\kappa\mu^{n-\frac12}$.
In this sense, we say that the energy dissipation law
in Theorem \ref{thm: CN energy law TFAC}
is asymptotically compatible in the limit $\alpha\rightarrow 1.$

\begin{remark}
 Other than the energy dissipation law \eqref{cont: TFAC energy law}, the TFAC model \eqref{cont: TFAC model} also admits the maximum bound principle \cite{TangYuZhou:2019}, that is, the solution is uniformly bounded by 1 if the initial
 and boundary data are uniformly bounded by 1. It is known that the numerical schemes \cite{LiaoTangZhou:2021,LiaoTangZhou:2020,LiaoZhuWang:2022NMTMA} based on 
 the nonuniform L1, L2-1$_\sigma$ and L1$_R$ formulas always preserve the maximum bound principle. 
 It is very interesting whether the nonuniform L1$^+$ time-stepping scheme \eqref{scheme: CN TFAC model} maintains the maximum  bound principle. This issue would be also challenging due to the lack of uniform monotonicity of the L1$^+$ kernels \eqref{def: L1+ formula} and remains open to us up to now. 
 \end{remark}

\section{Application to time-fractional Klein-Gordon model}
\setcounter{equation}{0}

Nonlinear integro-differential (fractional wave) equations play an important role
for describing anomalous diffusion processes and
wave propagation in viscoelastic materials \cite{Brunner2004book}.
We consider the following Klein-Gordon-type fractional wave equation \cite{AdolfssonEnelundLarsson:2003,GolmankhanehGolmankhanehBaleanu:2011,LyuVong:2022JSC}
with the fractional order $\beta\in(0,1)$,
\begin{align}\label{cont: TFKG model}
	\partial_tU+{\cal I}_t^{\beta}\zeta=0\quad\text{with
	$\zeta:=\frac{\delta E}{\delta U}=f(U)-\epsilon^2\Delta U$},
\end{align}
where $f(U)=F'(U)$ and
the associated kinetic energy $E[U]$ is defined as follows,
\begin{align}\label{cont: KG free energy}
	E[U]:=\int_{\Omega}\braB{\frac{\epsilon^2}{2}\abs{\nabla U}^2
		+F(U)}\zd\mathbf{x}\quad\text{with\; $F(U):=\frac{1}{4}\bra{U^2-1}^2$.}
\end{align}
This model \eqref{cont: TFKG model} is intermediate between the Allen-Cahn-type diffusion equation ($\beta= 0$)
and the Klein-Gordon-type wave equation ($\beta=1$),
and it can be termed as a nonlinear fractional PDE with the Caputo time derivative of order $\alpha=1+\beta\in(1,2)$.

Typically, in the limit $\beta\rightarrow 1$, the model \eqref{cont: TFKG model}
recovers the classical Klein-Gordon  equation $\partial_t^2U=\epsilon^2\Delta U-f(U)$,
which is a relativistic wave equation and describes the spin-zero particles in quantum field \cite{Greiner1994}.
As well-known, it admits the energy conservation law \cite{LiVu-Quoc:1995}
\begin{align}\label{KG-energy conservation}
	\frac{\zd {\cal E}}{\zd t}=0,
\end{align}
where the Hamiltonian energy ${\cal E}$ is defined by
\begin{align}\label{KG-energy}
	{\cal E}[U]:=E[U]+\frac12\|\partial_tU\|^2.
\end{align}
Therefore, it is natural to ask whether the time-fractional Klein-Gordon equation
\eqref{cont: TFKG model}  also maintains a similar energy law, and
whether the second-order time-stepping scheme based on
integral averaged formula \eqref{def: IAF formula} can also maintain
the corresponding energy law at the discrete time levels.
These problems seem also very important for the long-time numerical simulation
of certain integro-differential models containing nonlinear historical memory
terms (inside the time integral), cf. \cite{Brunner2004book}. To the best of our knowledge,
there were seldom related studies on the energy dissipation law
of nonlinear integro-differential models.

\subsection{Continuous energy dissipation law}

The first aim of this section is to define a new variational energy dissipation law
of \eqref{cont: TFKG model}. We consider the spatial domain $\Omega\subseteq \mathbb{R}^2$
with the periodic boundary conditions of $U$.
By the model \eqref{cont: TFKG model},
the time derivative of free energy
\begin{align}\label{cont: TFKG derivative-energy}
	\frac{\zd }{\zd t}E[U]=\myinnerb{\tfrac{\delta E}{\delta U},\partial_tU}
	=-\myinnerb{\zeta,{\cal I}_t^{\beta}\zeta}.
\end{align}
Taking $w:=\zeta$ in Lemma \ref{lem: countinous DGS equality}
with
$v=\mathcal{I}_t^{\beta}\zeta=-\partial_tU,$
one gets the following equality
\begin{align*}
	\myinnerb{\zeta,{\cal I}_t^{\beta}\zeta}
	=&\,\frac{1}2\,{}^R\!\partial_t^{\beta}\mynormb{\partial_tU}^2
	+\frac{1}{2}\int_0^t\frac{\partial}{\partial\xi}\bra{\frac{1}{\omega_{1-\beta}(t-\xi)}}
	\mynormB{\int_0^{\xi}\omega_{1-\beta}(t-s)v'(s)\zd{s}}^2\zd \xi.
\end{align*}
Thus by using \eqref{cont: TFKG derivative-energy} we get an energy dissipation law
\begin{align}\label{cont: TFKG energy law}
	\frac{\zd \mathcal{E}_\beta}{\zd t}+\frac{1}{2}\int_0^t\frac{\partial}{\partial\xi}\bra{\frac{1}{\omega_{1-\beta}(t-\xi)}}
	\mynormB{\int_0^{\xi}\omega_{1-\beta}(t-s)v'(s)\zd{s}}^2\zd \xi=0,
\end{align}
where the nonlocal energy
\begin{align}\label{cont: TFGK modified energy}
	\mathcal{E}_\beta[U]:=E[U]
	+\frac{1}{2}\mathcal{I}_{t}^{1-\beta}\mynormb{\partial_tU}^2\quad\text{for $t>0$.}
\end{align}
Obviously, as the fractional order $\beta\rightarrow1$, one has $\mathcal{E}_\beta[U]\rightarrow \mathcal{E}[U]$, see \eqref{KG-energy}.
Also, we can show that the second term in \eqref{cont: TFKG energy law}
vanishes as $\beta\rightarrow1$ (see the proof of Lemma \ref{lem: countinous DGS equality},
the difference term $J[v]$ naturally vanishes as the fractional
index $\beta\rightarrow1$). The energy dissipation law \eqref{cont: TFKG energy law}
degrades into \eqref{KG-energy conservation},
which is just the energy conservation law of the classical Klein-Gordon model.
In this sense, we say that both the nonlocal energy $\mathcal{E}_\beta[U]$
and the energy dissipation law \eqref{cont: TFKG energy law}  are asymptotically compatible
in the fractional order limit $\beta\rightarrow1$.

\subsection{Crank-Nicolson scheme}

Let $u^n\approx U(t_n)$ be the numerical solution.
By applying the integral averaged formula \eqref{def: IAF formula}
and the nonlinear approximation \eqref{def: nonlinear approximation},
we have the  following Crank-Nicolson scheme
\begin{align}\label{scheme: CN TFKG model}
	\partial_{\tau}u^{n}
	+(\mathcal{I}_{\tau}^{\beta}\zeta)^{n-\frac12}=0\quad\text{with}\quad
	\zeta^{n-\frac12}=f\brat{u}^{n-\frac12}-\epsilon^2\Delta u^{n-\frac12}
	\quad\text{for $n\ge 1$.}
\end{align}

Recalling the discrete kernels $a^{(\beta,n)}_{n-k}$ defined in \eqref{def: IAF formula kernel}
with the notation
$\mathcal{L}^{n-1}:=\sum_{k=1}^{n-1}a^{(\beta,n)}_{n-k}\tau_k\zeta^{k-\frac12}$,
we consider the following discrete functional $\mathcal{G}[z]$,
\begin{align*}
	\mathcal{G}[z]:=&\,\frac{1}{2a^{(\beta,n)}_{0}\tau_n^2}\mynormb{z-u^{n-1}}^2
	+\frac{1}{a^{(\beta,n)}_{0}\tau_n}\myinnerb{\mathcal{L}^{n-1},z}
	+\frac{\epsilon^2}{4}\mynormb{\nabla \brab{z+u^{n-1}}}^2\\
	&\,+\frac{1}{16}\mynormb{z}_{L^4}^4
	+\frac{1}{12}\myinnerb{u^{n-1},z^3}
	+\frac{1}{8}\myinnerb{\brat{u^{n-1}}^2,z^2}
	+\frac{1}{4}\myinnerb{\brat{u^{n-1}}^3,z}
	-\frac{1}4\mynormb{z+u^{n-1}}^2.
\end{align*}
The solution $u^n$ of nonlinear equation \eqref{scheme: CN TFKG model} at the time level $t_n$
is equivalent to the minimum of $\mathcal{G}[z]$ if and only if it is strictly convex and coercive.
With the necessary condition $a^{(\beta,n)}_{0}\tau_n^2\le 2$ for the convexity,
one can follow the proof of \cite[Theorem 2.1]{LiaoZhuWang:2022NMTMA} to prove the following result.

\begin{lemma}
	If $\tau_n\le \sqrt[1+\beta]{2\Gamma(2+\beta)}$,
	the Crank-Nicolson scheme \eqref{scheme: CN TFKG model} is uniquely solvable.
\end{lemma}

\subsection{Discrete energy dissipation law}

By using Theorem \ref{thm: DGS-IAF kernels},
we establish  a discrete energy  law for
the Crank-Nicolson scheme \eqref{scheme: CN TFKG model}.
With the original energy $E\kbra{u^{n}}$
in \eqref{cont: KG free energy}, we define the following discrete analogue of \eqref{cont: TFGK modified energy},
\begin{align}\label{def: variational energy TFKG}
	\mathcal{E}_\beta\kbrat{u^{n}}:=E\kbrat{u^{n}}
	+\frac{1}{2}\sum_{j=1}^n\mathsf{p}_{n-j}^{(\beta,n)}\mynormb{v^j}^2
	\quad\text{with $v^j:=\sum_{\ell=1}^{j}\mathsf{a}_{j-\ell}^{(\beta,j)}\tau_{\ell}\zeta^{\ell-\frac12}$.}
\end{align}

\begin{theorem}\label{thm: CN energy law TFKG}
	Under the step-ratio constraint \eqref{ieq: condition}, the 
	 variable-step Crank-Nicolson scheme \eqref{scheme: CN TFKG model}
	is energy stable in the sense that it
	preserves a discrete energy dissipation law
	\begin{align*}
		\partial_{\tau}\mathcal{E}_\beta\kbrat{u^n}
		+\frac{1}{2\tau_n}\sum_{j=1}^{n-1}\braB{\frac{1}{\Rcc{\mathsf{p}}_{n-j}^{(\beta,n)}}
			-\frac{1}{\Rcc{\mathsf{p}}_{n-j-1}^{(\beta,n)}}}
		\mynormB{\sum_{k=1}^{j}\Rcc{\mathsf{p}}_{n-k}^{(\beta,n)}\diff v^k}^2=0.
	\end{align*}
\end{theorem}
\begin{proof}
	Taking the inner product of  \eqref{scheme: CN TFKG model} by  $\tau_n\zeta^{n-\frac12}$, one gets
	\begin{align*}
		\myinnerb{\diff u^{n},\zeta^{n-\frac12}}
		+\myinnerb{(\mathcal{I}_{\tau}^{\beta}\zeta)^{n-\frac12},\tau_n\zeta^{n-\frac12}}=0.
	\end{align*}
	With the help of the equality \eqref{eq: property nonlinear approximation}, it is easy to check that
	$$\myinnerb{\diff u^{n},\zeta^{n-\frac12}}=E\kbrat{u^n}-E\kbrat{u^{n-1}},$$
	and we have
	\begin{align}\label{thmProof: KG energy law inner}
	E\kbrat{u^n}-E\kbrat{u^{n-1}}
	+\myinnerb{(\mathcal{I}_{\tau}^{\beta}\zeta)^{n-\frac12},\tau_n\zeta^{n-\frac12}}=0.
	\end{align}
	Taking $w_n:=\tau_n\zeta^{n-\frac12}$ in \eqref{eq: DGS equality}, we have
	\begin{align*}
		\myinnerb{(\mathcal{I}_{\tau}^{\beta}\zeta)^{n-\frac12},\tau_n\zeta^{n-\frac12}}=&\,
		\frac{1}{2}\sum_{j=1}^{n}\mathsf{p}_{n-j}^{(\beta,n)}\mynormb{v^j}^2
		-\frac{1}{2}\sum_{j=1}^{n-1}\mathsf{p}_{n-1-j}^{(\beta,n-1)}\mynormb{v^j}^2\nonumber\\
		&\,+\frac{1}{2}\sum_{j=1}^{n-1}\braB{\frac{1}{\Rcc{\mathsf{p}}_{n-j}^{(\beta,n)}}
			-\frac{1}{\Rcc{\mathsf{p}}_{n-j-1}^{(\beta,n)}}}
		\mynormB{\sum_{k=1}^{j}\Rcc{\mathsf{p}}_{n-k}^{(\beta,n)}\diff v^k}^2,
	\end{align*}
	where the sequence $v^j:=\sum_{\ell=1}^{j}\mathsf{a}_{j-\ell}^{(\beta,j)}\tau_{\ell}\zeta^{\ell-\frac12}$.
	Inserting it into \eqref{thmProof: KG energy law inner}
	yields the discrete energy dissipation law immediately.
	This completes the proof.
\end{proof}

As the fractional order $\beta\rightarrow1$,
	the definition \eqref{def: IAF formula kernel} gives the values
	$a_{0}^{(1,n)}=1/2$ and $a_{n-k}^{(1,n)}=1$ for $1\le k\le n-1$.
	The  time-stepping scheme \eqref{scheme: CN TFKG model}
	degrades into the Crank-Nicolson type scheme
	\begin{align}\label{scheme: CN KG model}
		\partial_{\tau} u^{n}
		+\frac{\tau_n}{2}\zeta^{n-\frac12}+\sum_{j=1}^{n-1}\tau_j\zeta^{j-\frac12}=0\quad\text{with}\quad
		\zeta^{n-\frac12}=f\brat{u}^{n-\frac12}-\epsilon^2\Delta u^{n-\frac12}
		\quad\text{for $n\ge 1$.}
	\end{align}
	One can check that this numerical scheme is uniquely solvable if $\tau_n\le 2$.
	This numerical scheme \eqref{scheme: CN KG model} 
	can be formulated into $\partial_{\tau} u^{n}+w^{n-\frac12}=0$ by introducing $w^n:=\sum_{k=1}^n\tau_k\zeta^{k-\frac12}$.
	With the fact $w^{n}-w^{n-1}=\tau_n\zeta^{n-\frac12}$, it is easy to establish 
	a discrete energy conservation law
	\[
	E\kbrat{u^n}
	+\frac1{2}\mynormb{w^n}^2=E\kbrat{u^{n-1}}
	+\frac1{2}\mynormb{w^{n-1}}^2
	\quad \text{for  $n\ge1$.}
	\]
	The definition \eqref{def: modified kernels} shows that the modified kernels
	$\mathsf{a}_{n-k}^{(1,n)}=1$ for $1\le k\le n$.
	The associated DOC kernels $\theta_{0}^{(1,n)}=1$,
	$\theta_{1}^{(1,n)}=-1$ and $\theta_{n-k}^{(1,n)}=0$ for $1\le k\le n-2$.
	Then the corresponding DCC  and RCC kernels read
	\begin{align*}
		&\mathsf{p}_{0}^{(1,n)}=1\quad \text{and}\quad\mathsf{p}_{n-k}^{(1,n)}=0
		\quad \text{for $1\le k\le n-1$},\\
		&\Rcc{\mathsf{p}}_{0}^{(1,n)}=1\quad \text{and}\quad\Rcc{\mathsf{p}}_{n-k}^{(1,n)}=0
		\quad \text{for $1\le k\le n-1$}.
	\end{align*}
	In this case, the term $\sum_{k=1}^{j}\Rcc{\mathsf{p}}_{n-k}^{(1,n)}(v_k-v_{k-1})$
	in the DGS \eqref{eq: DGS equality} naturally vanishes for $1\le j\le n-1$.
	Then  the DGS equality \eqref{eq: DGS equality} with $w_k:=\tau_k\zeta^{k-\frac12}$ degrades into
	\begin{align*}
		2w_n\sum_{j=1}^na^{(1,n)}_{n-j}w_j=v_n^2
		-v_{n-1}^2\quad\text{for $v_n:=\sum_{k=1}^n\tau_{k}\zeta^{k-\frac12}$ and $n\ge1$.}
	\end{align*}
	That is, the DGS equality \eqref{eq: DGS equality} is asymptotically compatible
	in the limit $\beta\rightarrow 1.$
	Obviously, the above discrete energy \eqref{def: variational energy TFKG} degrades into
	\begin{align*}
		\mathcal{E}_\beta\kbrat{u^{n}} \quad\longrightarrow\quad E\kbrat{u^{n}}
		+\frac1{2}\mynormb{v^n}^2\quad\text{as $\beta\rightarrow 1$;}
	\end{align*}
	and the discrete energy dissipation law in Theorem \ref{thm: CN energy law TFKG} degrades into
	\[
	\partial_\tau \braB{E\kbra{u^n}+\frac1{2}\mynormb{v^n}^2}=0\quad \text{for  $n\ge1$,}
	\]
	which is just the energy conservation law of \eqref{scheme: CN KG model}.
	Thus both the discrete energy \eqref{def: variational energy TFKG}
	and the energy dissipation law	in Theorem \ref{thm: CN energy law TFKG}
	are asymptotically compatible in the limit $\beta\rightarrow 1.$

\section{Numerical experiments}
\setcounter{equation}{0}

The Fourier pseudo-spectral method is employed for the spatial discretization. 
The spatial domain $\Omega=(0,2\pi)^2$ is discretized by using
$32\times32$ uniform grids.
The resulting nonlinear system at each time level
is solved by using a simple fixed-point iteration with the termination error $10^{-12}$.
The sum-of-exponentials technique \cite{JiangZhangQianZhang:2017,LiaoTangZhou:2020}
with an absolute tolerance error $\varepsilon=10^{-12}$ and cut-off time $\Delta{t}=\tau_1$
is always adopted to reduce the computational cost and storage. As done in \cite{LiaoMcLeanZhang:2021alikhanov,LiaoTangZhou:2020}, we always compute the discrete coefficients 
	\eqref{def: IAF formula kernel} and \eqref{def: L1+ formula kernel} 
	with adaptive Gauss-Kronrod quadrature  to avoid the roundoff error problem.

\begin{figure}[htb!]
	\centering
	\subfigure[The regularity parameter $\sigma=0.4$]{
		\includegraphics[width=2.0in]{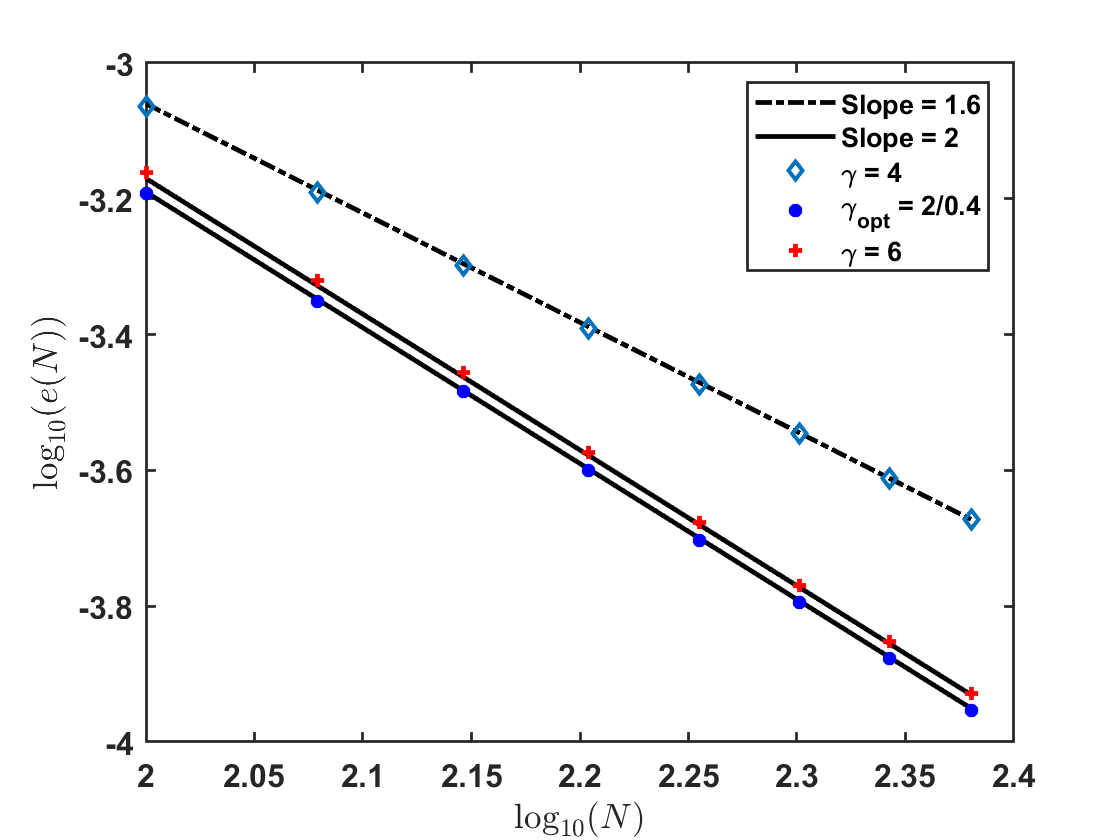}}
	\subfigure[The regularity parameter $\sigma=0.8$]{
		\includegraphics[width=2.0in]{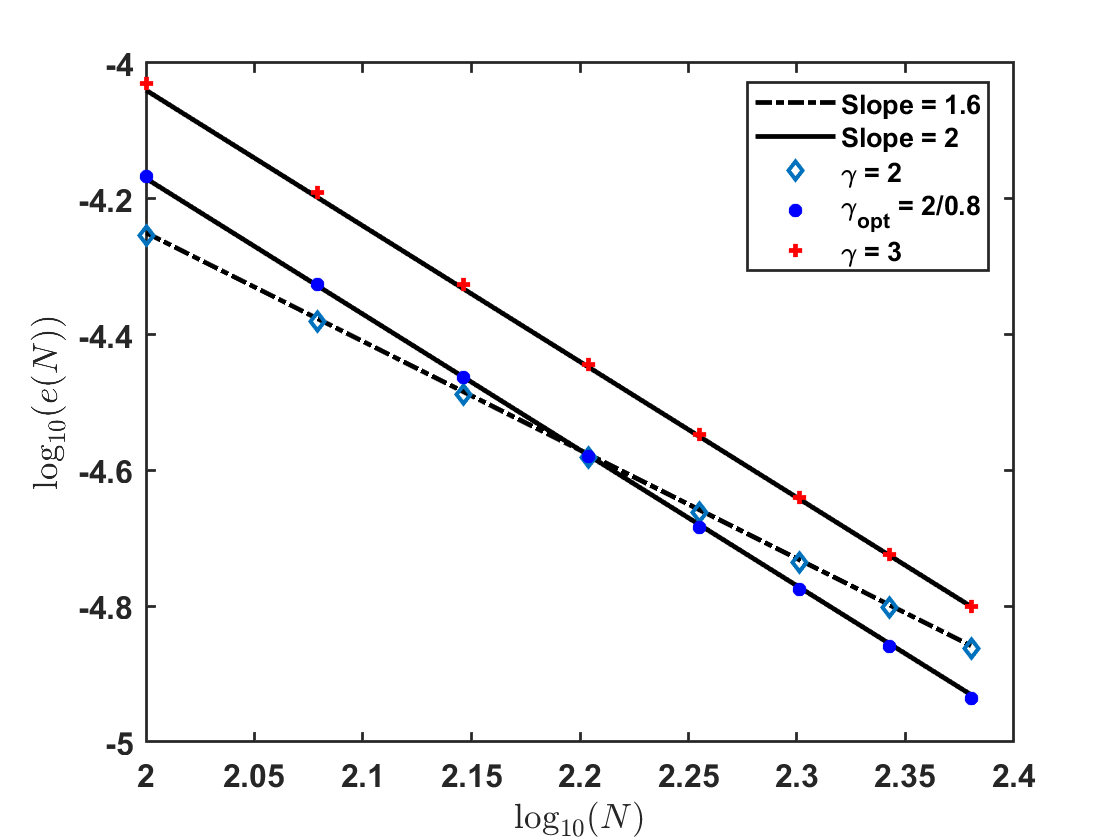}}
	\caption{The log-log plot of convergence with different regularity parameters $\sigma$. }
	\label{singu:L1 AC alpha 08}
\end{figure}

\begin{example}
	We solve the forced model $\partial_{t}^{\alpha}\Phi=-\kappa(\mu-g)$ 
	by adding a forcing term $g=g(\mathbf{x},t)$ to 
	the time-fractional Allen-Cahn model \eqref{cont: TFAC model}
	with $\kappa=1$ and $\epsilon^2=0.5$
	for $\mathbf{x}\in(0,2\pi)^{2}$ and $0<t\le T=1$
	such that the exact solution $\Phi=\omega_{1+\sigma}(t)\sin x\sin y$,
	where $\sigma\in(0,1)$ is a regularity parameter.
\end{example}

We take the graded time mesh $t_k=T_0(k/N_0)^\gamma$ for $0\le k\le N_0$
in the interval $[0, T_0]$, where
$T_0=\min\{1/\gamma,T\}$ and $N_0=\lceil \frac{N}{T+1-\gamma^{-1}}\rceil$.
In the remainder interval $[T_{0},T]$,
the random time meshes
$\tau_{N_{0}+k}:=(T-T_{0})s_{k}/S_1$ for $1\le k\le N_1$
are used
by setting $N_1:=N-N_0$ and $S_1=\sum_{k=1}^{N_1}s_k$,
where $s_k\in(0,1)$ are random numbers.
By taking the fractional order $\alpha=0.8$ 
(the results for other fractional orders are similar and omitted here),
we run the Crank-Nicolson scheme \eqref{scheme: CN TFAC model} 
for different total numbers $N=100+20m$ $(1\le m\le 8)$
and record the $L^2$ norm error $e(N):=\max_{1\le{n}\le{N}}\mynorm{\Phi^n-\phi^n}$
in each run. With three different grading parameters $\gamma$ with $\gamma_{\text{opt}}:=2/\sigma$, 
Figure \ref{singu:L1 AC alpha 08} depicts 
the experimental order of convergence in the log-log plot for two different regularity parameters 
$\sigma=0.4$ and $0.8$. We observe that the time accuracy is of order $O\brat{\tau^{\gamma\sigma}}$
when the graded parameter $\gamma<\gamma_{\text{opt}}$;
while the second-order accuracy $O\brat{\tau^{2}}$ can be achieved when the graded parameter
$\gamma\ge\gamma_{\text{opt}}$.

\begin{example}\label{example:coarsen dynamic}
	We simulate the coarsening dynamics of the time-fractional Allen-Cahn model \eqref{cont: TFAC model}.
	The initial condition is taken as
	$\phi_0(\mathbf{x})=\text{rand}(\mathbf{x})$,
	where $\text{rand}(\mathbf{x})$ generates uniform random numbers
	between $-0.001$ to $0.001$.
	The mobility coefficient $\kappa=1$ and
	the interfacial thickness $\epsilon=0.01$.
\end{example}

\begin{figure}[htb!]
	\centering
	\subfigure[Original energy $E\kbra{\phi^n}$]{
		\includegraphics[width=2.0in]{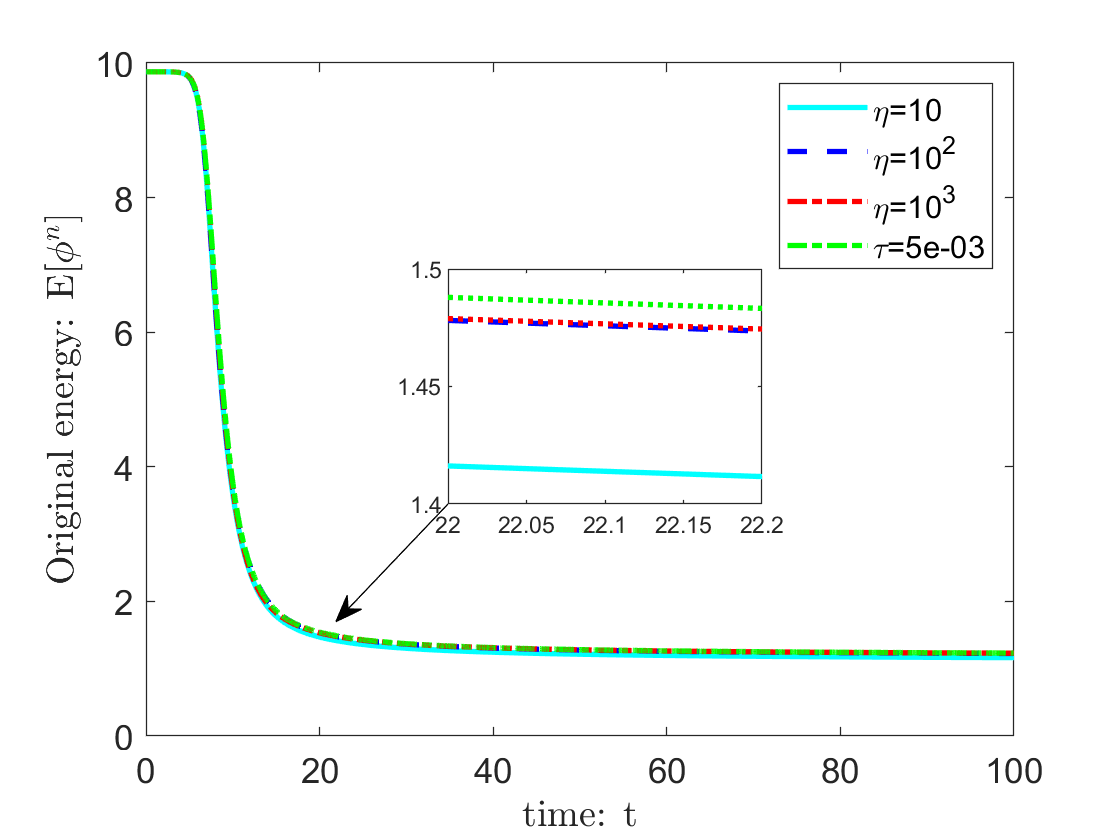}}
	\subfigure[Modified energy $E_\alpha\kbra{\phi^n}$]{
		\includegraphics[width=2.0in]{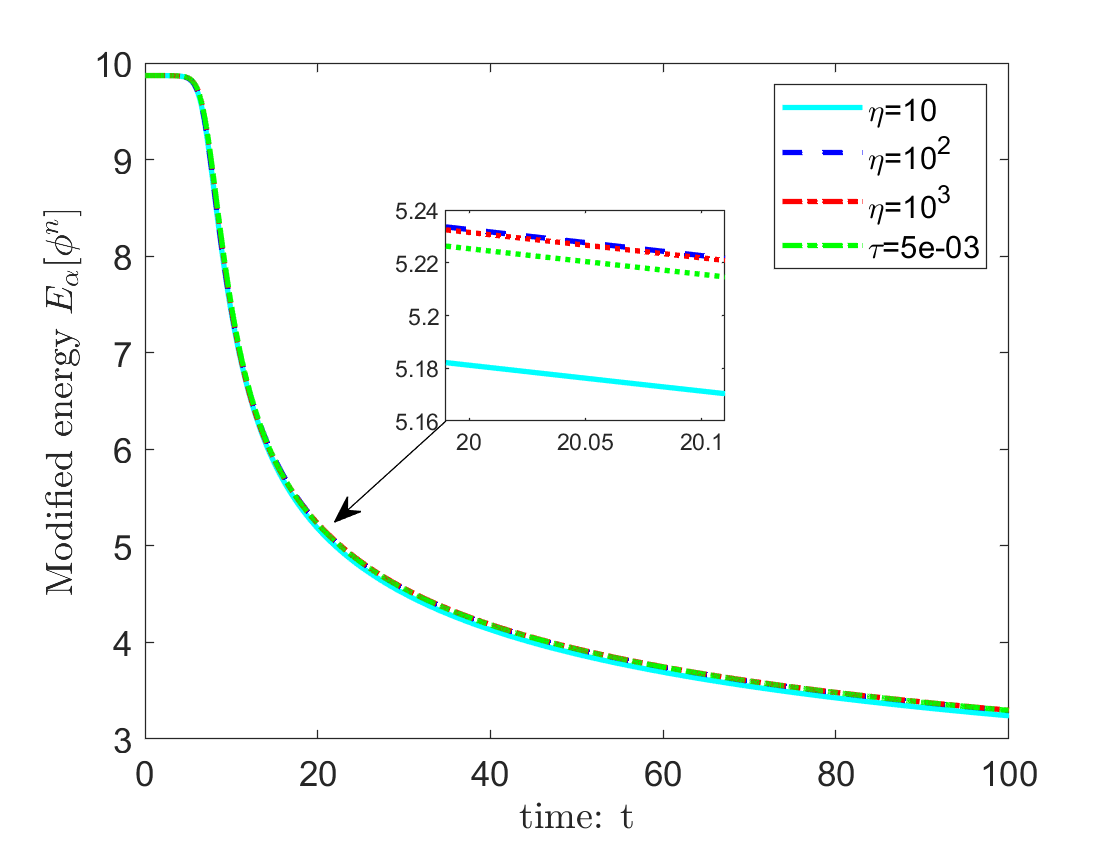}}
	\subfigure[Time steps $\tau_n$]{
		\includegraphics[width=2.0in]{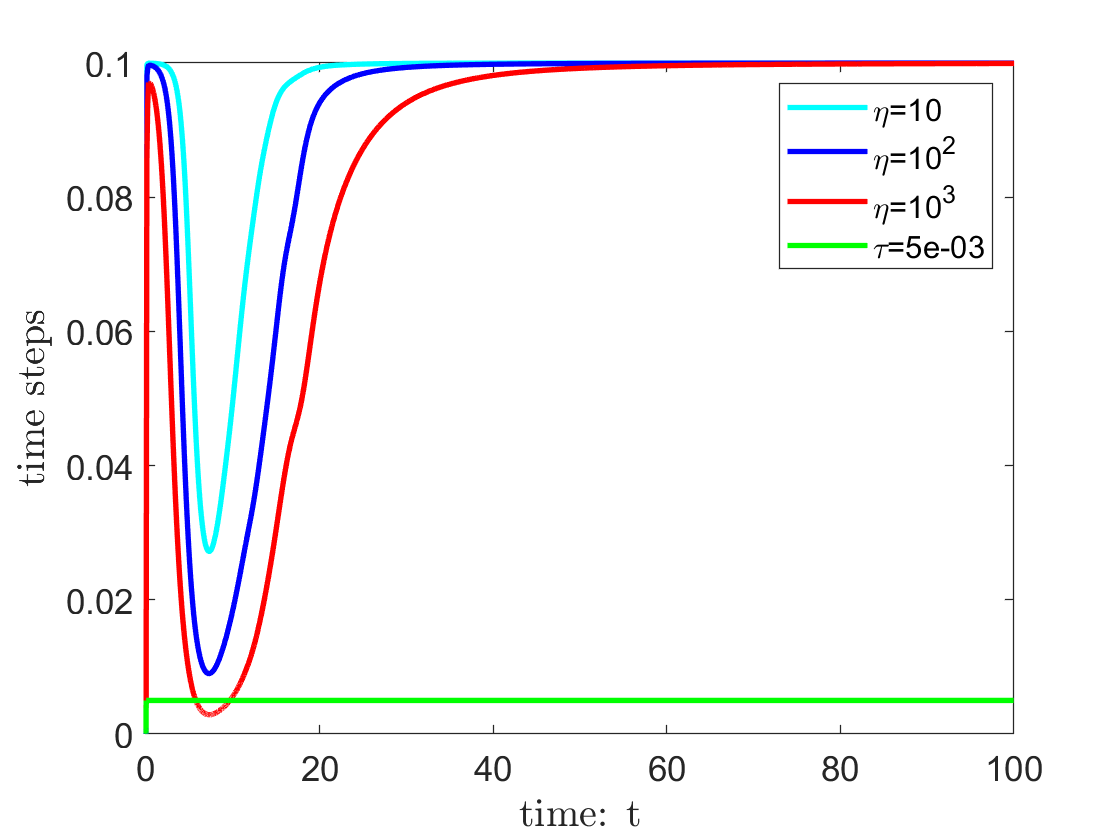}}
	\caption{Energy curves by uniform step and adaptive strategy with different parameters $\eta$.}
	\label{example:compar adap param}
\end{figure}

\begin{figure}[htb!]
	\centering
	\subfigure[Original energy $E\kbra{\phi^n}$]{
		\includegraphics[width=2.0in]{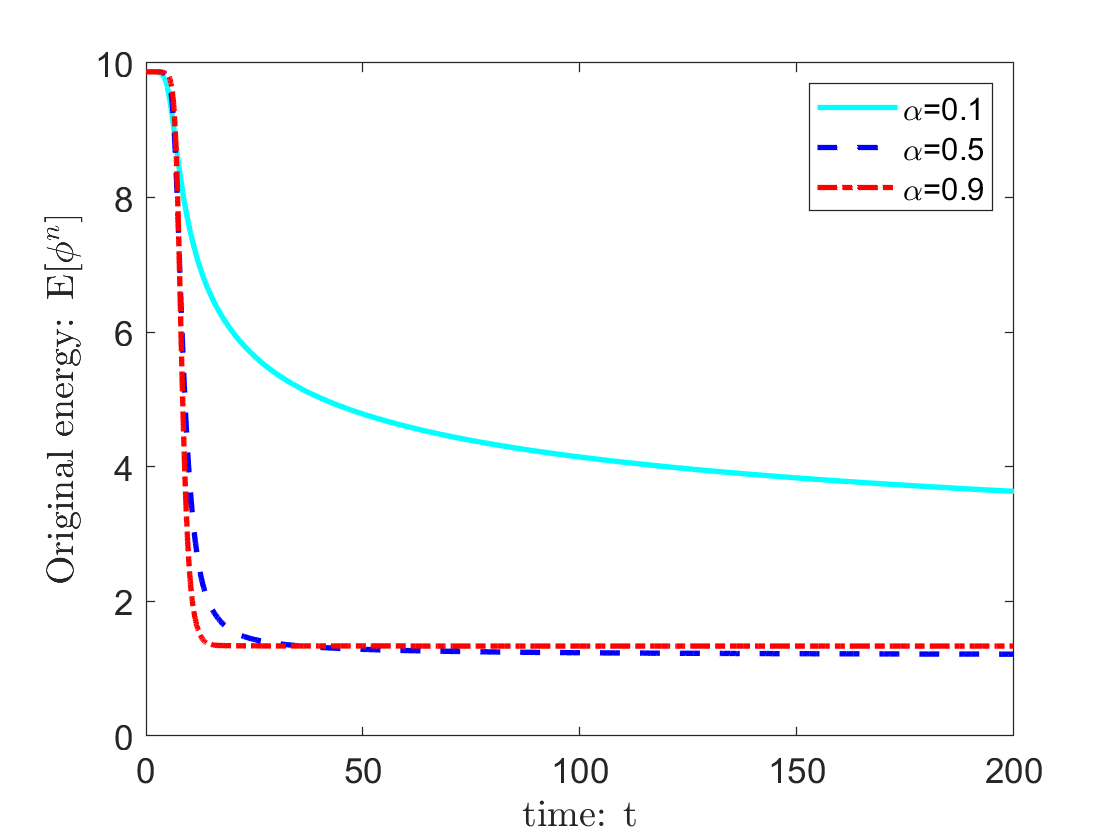}}
	\subfigure[Modified energy $E_\alpha\kbra{\phi^n}$]{
		\includegraphics[width=2.0in]{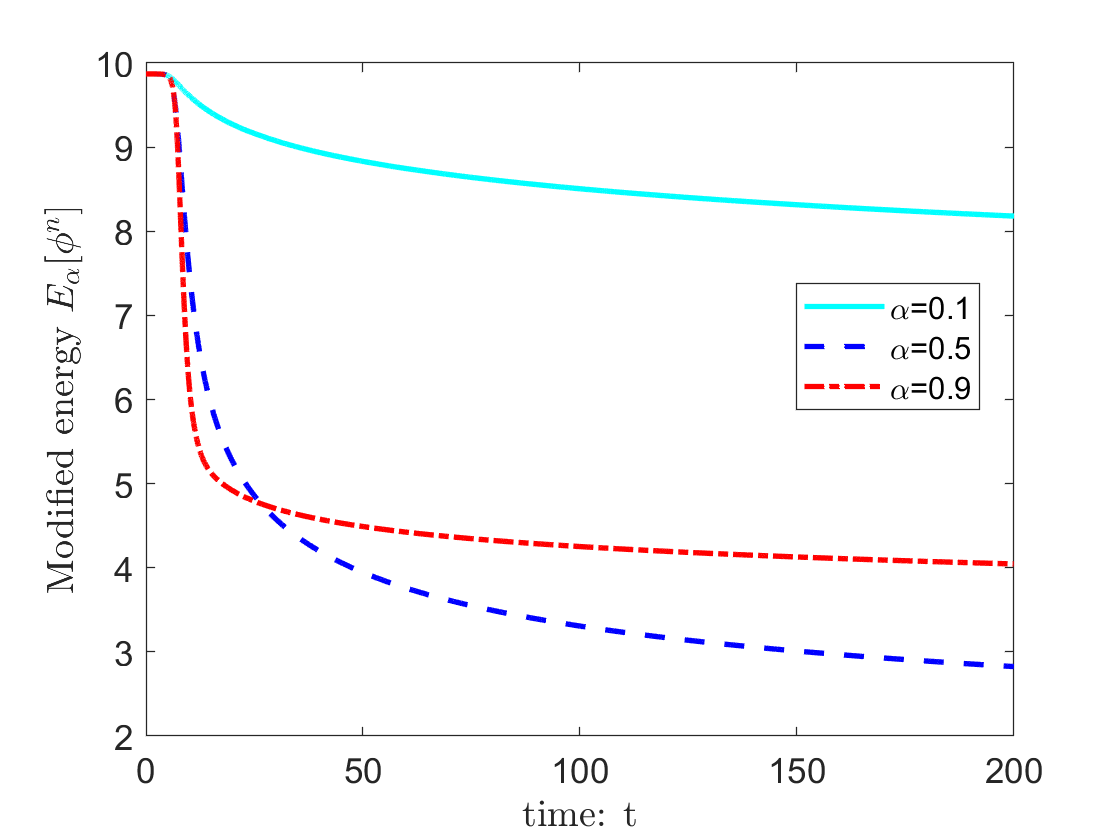}}
	\subfigure[Adaptive time steps $\tau_n$]{
		\includegraphics[width=2.0in]{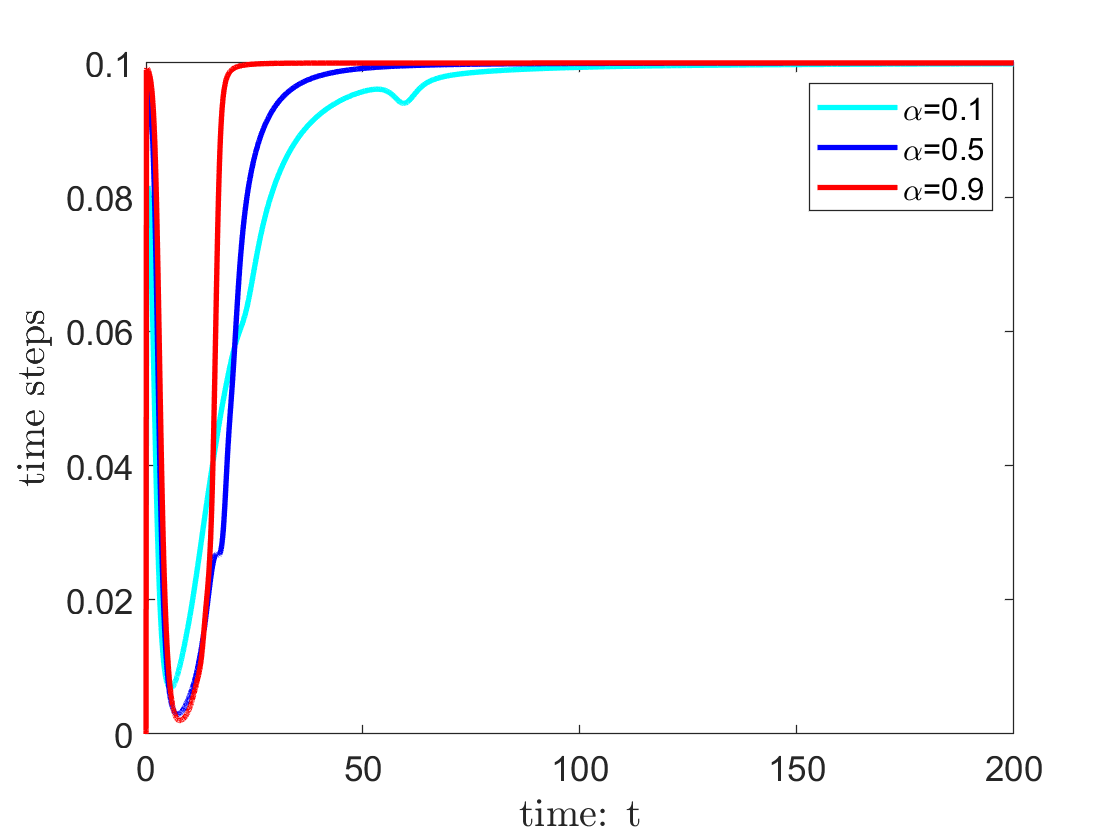}}
	\caption{Energy curves of  Example 2 with different fractional orders $\alpha$.}
	\label{example:FAC coarsen energy}
\end{figure}

The graded time mesh $t_k=T_0(k/N_0)^{\gamma}$ together with
the settings $\gamma=4,N_0=50$ and $T_0=0.01$ is applied
to resolve the initial singularity.
In the remainder interval $[T_0,T]$, we adjust the time-step sizes by 
the following time-stepping strategy \cite{JiZhuLiao:2022,LiaoZhuWang:2022NMTMA},
\begin{align}\label{algorithm:adaptive time stepping}
	\tau_{n+1}=\max\big\{\tau_{ada},r_{g}(r_n)\tau_n\big\}\quad\text{with}\quad 
	\tau_{ada}
	=\max\left\{\tau_{\min},
	\tau_{\max}/\Pi_{\phi}\right\}
\end{align}
where  $\tau_{\max}=0.1$ and $\tau_{\min}=10^{-3}$ are
the predetermined maximum and minimum time steps, respectively,  $\Pi_{\phi}:=\sqrt{1+\eta\mynormb{\partial_\tau\phi^n}^2}$ for 
a user parameter $\eta>0$, and $r_{g}$ is defined by \eqref{ieq: guess condition}. 
At first, we perform a comparative study by running
the scheme \eqref{scheme: CN TFAC model} with the fractional order $\alpha=0.5$.
A small uniform time step $\tau=5\times10^{-3}$ is used to compute the reference solution.
Figure \ref{example:compar adap param} plots the (original and modified) energy curves and the associated time-steps by
using the adaptive time-stepping strategy \eqref{algorithm:adaptive time stepping}  with three user parameters $\eta=10$, $10^2$ and $10^3$, respectively. It seems that the result of $\eta=10^3$ is well accordant with the reference solution.

We perform the numerical simulation by using
the adaptive time-stepping strategy \eqref{algorithm:adaptive time stepping}
with the user parameter $\eta=10^3$ until time $T=200$.
The curves of original energy $E[\phi^n]$ in \eqref{cont: classical free energy} 
and the modified energy $E_{\alpha}[\phi^n]$ in \eqref{cont: TFAC modified energy} together with 
the associated time steps during the coarsening dynamics are
depicted in Figure \ref{example:FAC coarsen energy}.
The numerical results are consistent with  those reported in \cite{LiaoTangZhou:2021,LiaoZhuWang:2022NMTMA,TangYuZhou:2019}
and our proposed method \eqref{scheme: CN TFAC model} can effectively
capture the multiple time scales in the long-time dynamical simulations.
We see that the value of fractional order $\alpha$ significantly affects
	the coarsening dynamics process, but it hardly affects the steady-state solution.

	\begin{figure}[htb!]
		\centering
		\subfigure[The regularity parameter $\sigma=0.4$]{
			\includegraphics[width=2.0in]{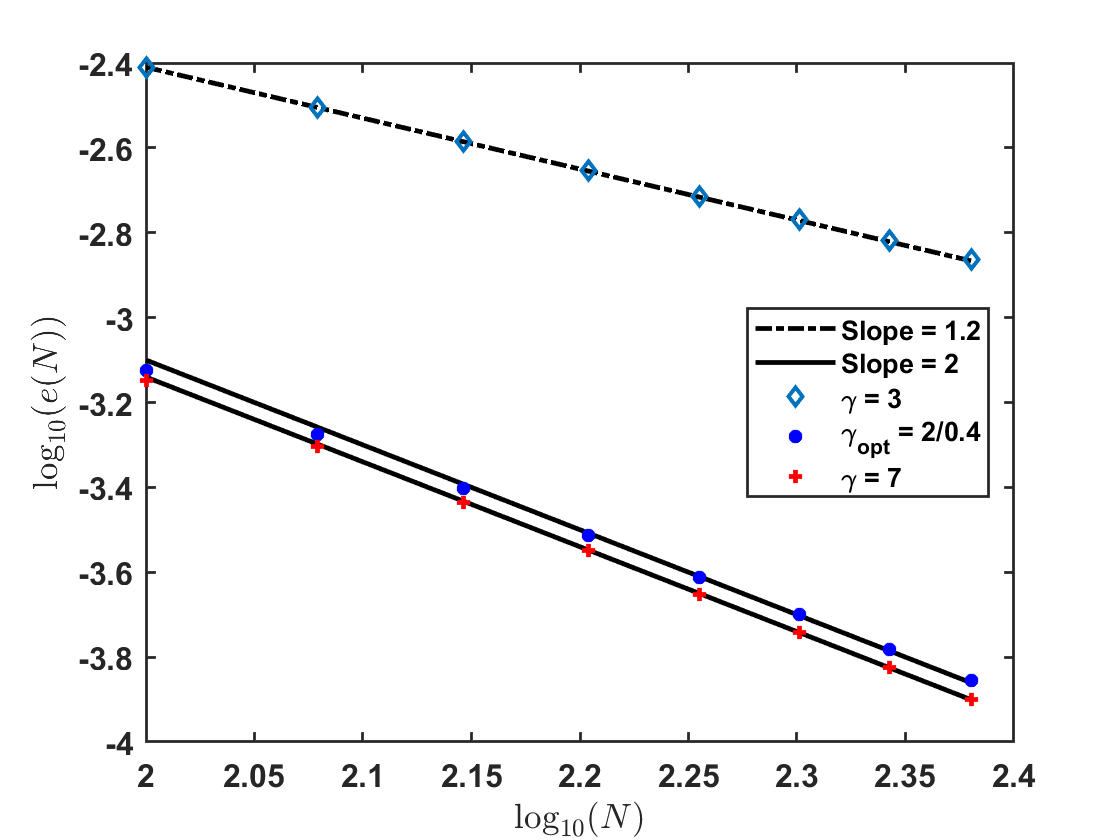}}
		\subfigure[The regularity parameter $\sigma=0.8$]{
			\includegraphics[width=2.0in]{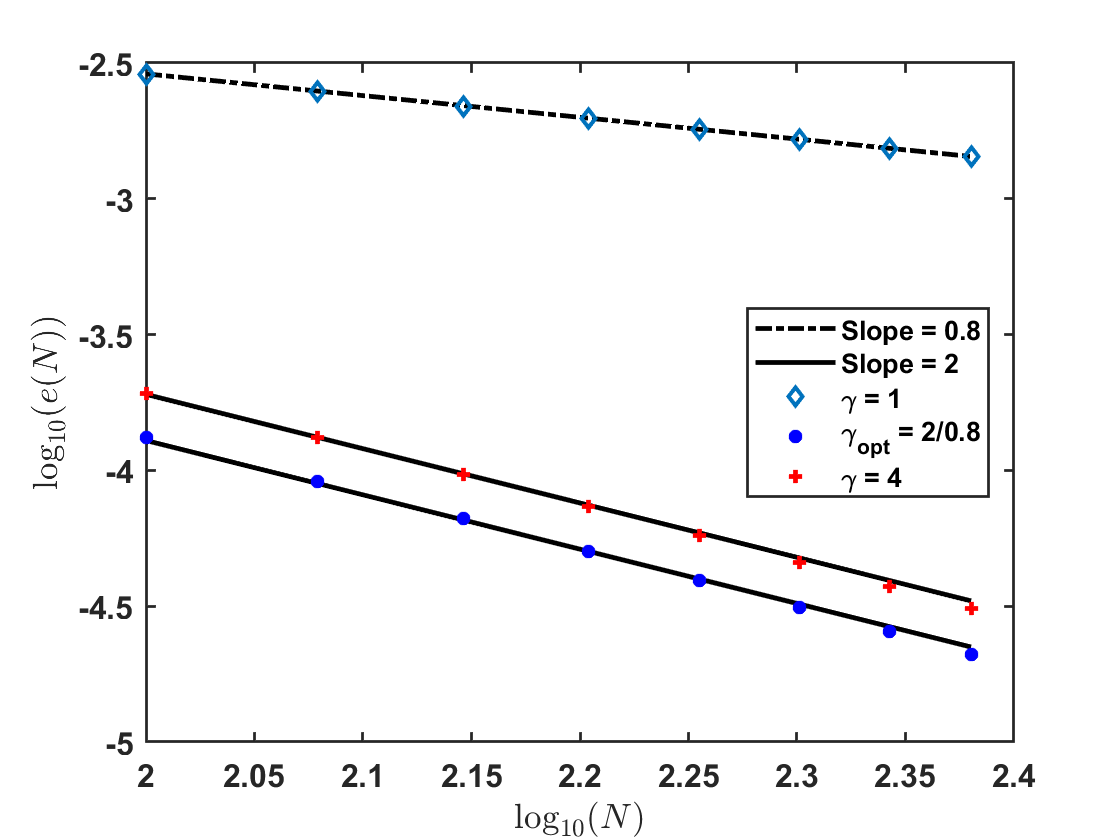}}
		\caption{The log-log plot of convergence with different regularity parameters $\sigma$. }
		\label{singu:L1 KG beta 05}
	\end{figure}
	
	\begin{example}
		We solve the forced model $\partial_tU+{\cal I}_t^{\beta}\zeta=g$ 
		by adding a forcing term $g=g(\mathbf{x},t)$ to 
		the time-fractional Klein-Gordon model \eqref{cont: TFKG model}
		with  $\epsilon^2=0.5$ for $\mathbf{x}\in(0,2\pi)^{2}$ and $0<t\le T=1$
		such that the exact solution $U=\omega_{1+\sigma}(t)\sin x\sin y$,
		where $\sigma\in(0,1)$ is a regularity parameter.
	\end{example}

\begin{figure}[htb!]
	\centering
	\subfigure[Modified energy $\mathcal{E}_\beta\kbra{u^n}$]{
		\includegraphics[width=2.0in]{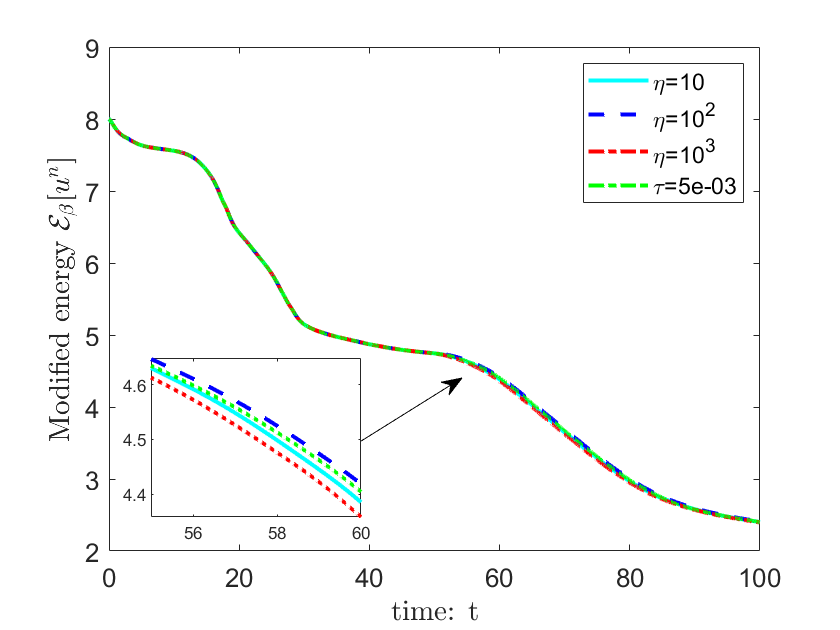}}
	\subfigure[Adaptive time steps $\tau_n$]{
		\includegraphics[width=2.0in]{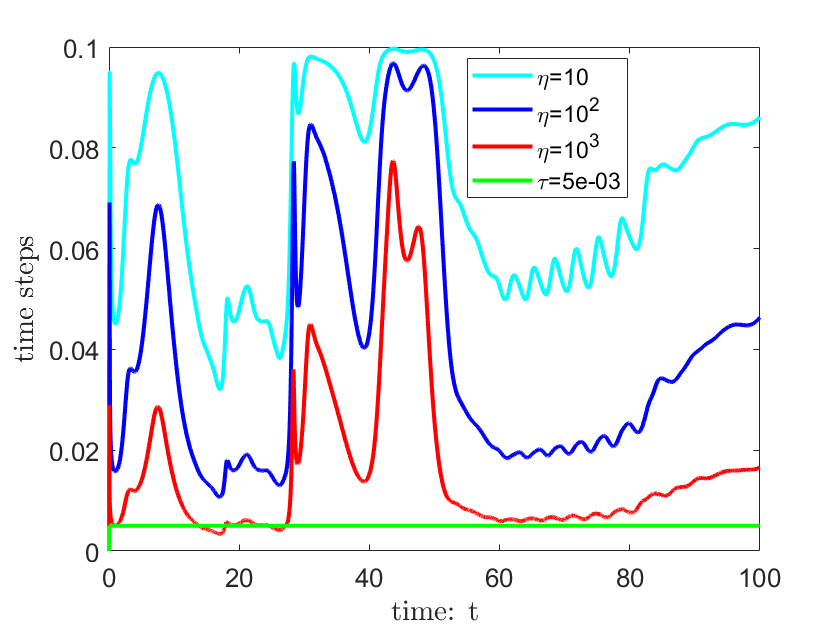}}
	\caption{Modified energy curves by adaptive time-steppings with different parameters $\eta$.}
	\label{example:FKG energy different parameter}
\end{figure}	
	\begin{figure}[htb!]
		\centering
		\subfigure[Modified energy $\mathcal{E}_\beta\kbra{u^n}$]{
			\includegraphics[width=2.0in]{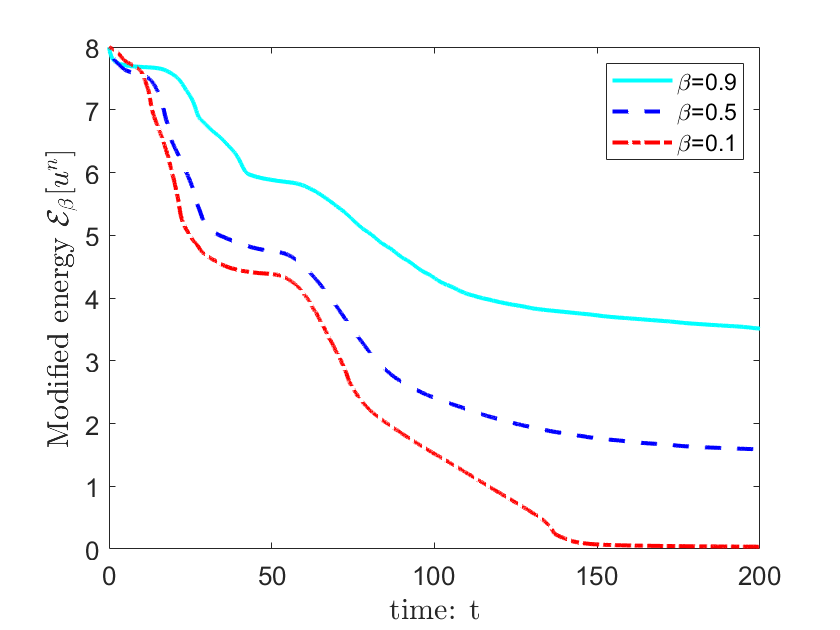}}
		\subfigure[Adaptive time steps $\tau_n$]{
			\includegraphics[width=2.0in]{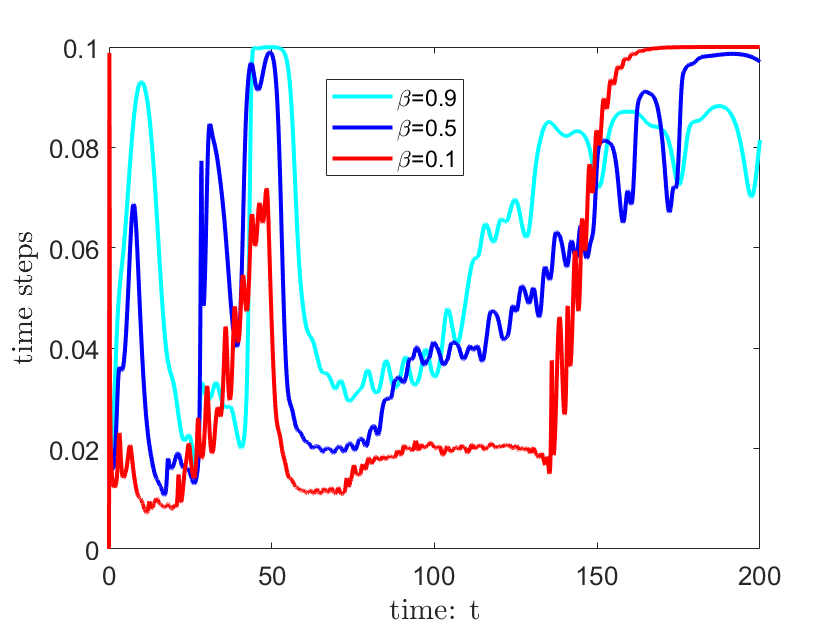}}
		\caption{Modified energy curves of  Example 4 with different fractional orders $\beta$.}
		\label{example:FKG energy different order}
	\end{figure}
	
	By taking the fractional index $\beta=0.5$,
	we run the Crank-Nicolson scheme \eqref{scheme: CN TFKG model} 
	for different total numbers $N=100+20m$ $(1\le m\le 8)$
	and record the $L^2$ norm error $e(N):=\max_{1\le{n}\le{N}}\mynorm{U^n-u^n}$
	in each run. With three different grading parameters $\gamma$ with $\gamma_{\text{opt}}:=2/\sigma$, 
	Figure \ref{singu:L1 KG beta 05} depicts 
	the experimental order of convergence in the log-log plot for two different regularity parameters 
	$\sigma=0.4$ and $0.8$. Again, we observe that the time accuracy is of order $O\brat{\tau^{\gamma\sigma}}$
	when the graded parameter $\gamma<\gamma_{\text{opt}}$;
	while the second-order accuracy $O\brat{\tau^{2}}$ can be achieved when the parameter
	$\gamma\ge\gamma_{\text{opt}}$.
	
	\begin{example}
		We simulate the time-fractional Klein-Gordon model \eqref{cont: TFKG model}
		with the parameter $\epsilon=0.1$ subject to the initial condition $u_0(\mathbf{x})=\cos3x\cos5y$.
	\end{example}
	
	Figure \ref{example:FKG energy different parameter} plots the curves 
	of the modified energy $\mathcal{E}_\beta[u^n]$ defined in \eqref{cont: TFGK modified energy} with the fractional index $\beta=0.5$ and three user parameters $\eta=10$, $10^2$ and $10^3$, respectively. It seems that the result of $\eta=100$ is well accordance with the reference solution generated with the small time-step $\tau=0.005$. The graded mesh $t_k=T_0(k/N_0)^{\gamma}$ with
	the settings $\gamma=4,N_0=50$ and $T_0=0.01$ is also applied
	to resolve the initial singularity. In the remainder interval $[T_0,T]$, we adopt 
	the adaptive time-stepping strategy \eqref{algorithm:adaptive time stepping}. 	
	We perform the numerical simulation until time $T=200$ with the user parameter $\eta=100$.
	The curves of the modified energy $\mathcal{E}_\beta[u^n]$ for different fractional orders  and 	
	the associated time steps are depicted in Figure \ref{example:FKG energy different order}. 
	It is seen that the adaptive time-stepping strategy can effectively
	capture the multiple time scales in long-time simulations.

\section*{Acknowledgements}
The authors would like to thank the editor and the anonymous referees for their valuable suggestions. 
They are helpful in improving the quality of the paper.

\appendix
\section{Some technical lemmas}
\setcounter{equation}{0}

\begin{lemma}\label{lem: app-hx}
	For the function $\rho$ in \eqref{def: app-varphi}, it holds that
	$$\rho(z)\le 2\kbrat{(z+1)^{\beta}-1}\quad\text{for $z\ge1$.}$$
\end{lemma}
\begin{proof}
	We define a function $h_1(z):=\rho(z)-2\kbra{(z+1)^{\beta}-1}$ with the first derivative
	\begin{align*}
		h_1'(z)=&\rho'(z)-2\beta(z+1)^{\beta-1}
		= \frac{1+\beta}{(z+1)^{1-\beta}}\kbraB{(1+z)-z(1+1/z)^{1-\beta}-\frac{2\beta}{1+\beta}}\\
		<&\,\frac{\beta(1-\beta)}{2z(z+1)^{1-\beta}}\bra{1+\beta-2z}<0\quad\text{for $z>1$,}
	\end{align*}
	where the following estimate due to the Taylor expansion has been used,
	\begin{align*}
		z(1+1/z)^{1-\beta}>z+1-\beta-\frac{(1-\beta)\beta}{2z}\quad\text{for $z>1$.}
	\end{align*}	
	Thus $h_1(z)\leq h_1(1)=0$ for $z\geq 1$, and the result is verified.
\end{proof}

\begin{lemma}\label{lem: app-varphi}
	For the function $\rho$ in \eqref{def: app-varphi}, it holds that
		$$\rho(z)< 2(2^\beta-1)z^\beta\quad\text{ for $0<z<1$.}$$	
\end{lemma}
\begin{proof}
	Let $h_2(z):=\rho(z)- 2(2^\beta-1)z^\beta$ with the first derivative
	$$h_2'(z)=(1+\beta)((z+1)^\beta-z^\beta)-2\beta(2^\beta-1)z^{\beta-1} \quad\text{for $z>0$}.$$
	To examine the property of $h_2'(z)$, we consider another auxiliary function $$h_3(z):=\ln\kbrab{(1+\beta)((z+1)^\beta-z^\beta)}-\ln\kbrab{2\beta(2^\beta-1)z^{\beta-1}} \quad\text{for $z>0$}. $$
	The first derivative 
	$$h_3'(z)=(\beta-1)\frac{\int_0^1(z+s)^{\beta-2}\zd s}{\int_0^1(z+s)^{\beta-1}\zd s}+\frac{1-\beta}{z}
	=\frac{1-\beta}{z}-\frac{1-\beta}{z+\xi}>0\quad\text{with}\quad 0<\xi<1,$$
	where the Cauchy differential mean-value theorem was applied in the second equality.
	We see that $h_3(z)$ is increasing with respect to $z$ for $z\in(0,1)$. Note that,  $$h_3(0)=\ln(1+\beta)-\lim_{z\rightarrow 0^+}\ln[2\beta(2^\beta-1)z^{\beta-1}]<0$$ 
	and $h_3(1)=\ln[(1+\beta)/(2\beta)]>0.$ 
	We know that $h_3(z)=0$ has a unique root $z_1\in(0,1)$ such that $h_3(z)<0$ for $0< z<z_1$ 
	and $h_3(z)>0$ for $z_1<z< 1$. 
	Thus, $h_2(z)$ is decreasing for $0< z<z_1$ and increasing for $z_1<z< 1$. Then we obtain 
	$h_2(z)< \max\{h_2(0),h_2(1)\}=0$ for $0< z< 1$ and complete the proof.
\end{proof}

\begin{lemma}\label{lem: app-Hxy}
	For the functions $\rho$ and $r_*$ defined in \eqref{def: app-varphi} and \eqref{ieq: condition}, respectively, 
	it holds that  
	\begin{align*}
		\frac{\rho(x)}{2x}< \frac{\rho(xy+y)-\rho(y)}{x\rho(y)} \quad\text{for $x\geq r_*(y)$ and $y>0$}. 
	\end{align*}
\end{lemma}
\begin{proof}Note that, $\rho'(z)=(1+\beta)\kbrat{(z+1)^{\beta}-z^{\beta}}>0$ and $$\rho''(z)=\beta(1+\beta)\kbrat{(z+1)^{\beta-1}-z^{\beta-1}}<0\quad\text{for $z>0$.}$$
	We consider an auxiliary function $h_4(x,y):=\rho(xy+y)$ with 
	$$\partial_xh_{4}=y\rho'(xy+y)>0\quad\text{and}\quad
	\partial_{xx}h_{4}=y^2\rho''(xy+y)<0\quad\text{for $x,y>0$}$$
	such that $h_4(x,y)$ is increasing and concave with respect to $x>0$. Then we get
	$$h_4(x,y)-h_4(0,y)>[h_4(1,y)-h_4(0,y)]x,$$
	 or
	$$\rho(xy+y)-\rho(y)> [\rho(2y)-\rho(y)]x\quad\text{for $0< x< 1$.}$$
	Since $r_*(y)<1$ for $y>0$ according to \eqref{ieq: r* condition},
	we apply Lemma \ref{lem: app-varphi} to derive that
	\begin{align*}
	\frac{\rho(x)}{2x}<(2^\beta-1)x^{\beta-1}\le \frac{\rho(2y)-\rho(y)}{\rho(y)}<\frac{\rho(xy+y)-\rho(y)}{x\rho(y)}
		\quad\text{for $r_*(y)\le x<1$}. 
	\end{align*}
The desired inequality is verified and the proof is completed. 	
\end{proof}


\begin{thebibliography}{99}
	
	\bibitem{AlsaediAhmadKirane:2015}
	A. Alsaedi, B. Ahmad and M. Kirane, Maximum principle for certain generalized time and space-fractional diffusion equations, Quart. Appl. Math., 73 (2015), pp. 163--175.
	
	\bibitem{AdolfssonEnelundLarsson:2003}
	K. Adolfsson, M. Enelund and S. Larsson, Adaptive discretization of an integro-differential equation with a weakly singular convolution kernel, Comput. Methods Appl. Mech. Engrg., 192 (2003), pp. 5285--5304.
	
	\bibitem{Brunner2004book}
	H. Brunner, Collocation Methods for Volterra Integral and Related Functional Equations,
	 Cambridge University Press, 2004.
	
	
	\bibitem{CuestaLubichPalencia:2006}
	E. Cuesta, C. Lubich and C. Palencia, Convolution quadrature time discretization of fractional diffusion-wave equations, Math. Comput., 75 (2006), pp. 673--696.
	
	
	\bibitem{CuestaPalencia:2003}
	E. Cuesta and C. Palencia, A numerical method for an integro-differential equation with memory in Banach spaces: qualitative properties, SIAM J. Numer. Anal., 41 (2003), pp. 1232--1241.
	
	\bibitem{GolmankhanehGolmankhanehBaleanu:2011}
	A. Golmankhaneh, A. Golmankhaneh and D. Baleanu, On nonlinear fractional Klein-Gordon equation, Signal Processing, 91 (2011), pp. 446--451.
	
	
	\bibitem{Greiner1994}
	W. Greiner, Relativistic Quantum Mechanics: Wave Equation, Springer, 1994.
    
	\bibitem{JiLiaoGongZhang:2020}
	B. Ji, H.-L. Liao, Y. Gong and L. Zhang, Adaptive second-order Crank--Nicolson time-stepping schemes for time-fractional molecular beam epitaxial growth models, SIAM J. Sci. Comput., 42 (2020), pp. B738--B760.
	
	\bibitem{JiZhuLiao:2022}
	B. Ji, X. Zhu and H.-L. Liao, Energy stability of variable-step L1-type schemes 
	for time-fractional Cahn-Hilliard model, Commun. Math. Sci., 2023, to appear.
	
	\bibitem{JiangZhangQianZhang:2017}
	S. Jiang, J. Zhang, Q. Zhang and Z. Zhang, Fast evaluation of the Caputo fractional
	derivative and its applications to fractional diffusion equations, Comm. Comput. Phys.,
	21 (2017), pp. 650--678.
	
	\bibitem{LiVu-Quoc:1995}
	S. Li and L. Vu-Quoc, Finite difference calculus invariant structure of a class of algorithms for the nonlinear Klein-Gordon equation, SIAM J. Numer. Anal., 32 (1995), pp. 1839--1875.
	
	
	
	
	\bibitem{LiaoMcLeanZhang:2019}
	H.-L. Liao, W. McLean and J. Zhang,  A discrete Gr\"{o}nwall inequality with applications to numerical schemes for subdiffusion problems, SIAM J. Numer. Anal., 57 (2019), pp. 218--237.
	
	\bibitem{LiaoMcLeanZhang:2021alikhanov}
	{ H.-L. Liao, W. McLean and J. Zhang},
	{ A second-order scheme with nonuniform time steps
		for a linear reaction-subdiffusion equation},
	Commun. Comput. Phys., 30 (2021), pp. 567--601.
	
	\bibitem{LiaoTangZhou:2021}
	H.-L. Liao, T. Tang and T. Zhou, An energy stable and maximum bound preserving scheme
	with variable time steps for time fractional Allen-Cahn equation, SIAM J. Sci. Comput., 43 (2021), pp. A3503--A3526.
	
	\bibitem{LiaoTangZhou:2020}
	H.-L. Liao, T. Tang and T. Zhou, A second-order and nonuniform time-stepping maximum-principle preserving scheme for time-fractional Allen-Cahn equations, J. Comput. Phys., 141 (2020), num. 109473, doi: 10.1016/j.jcp.2020.109473.
	
	\bibitem{LiaoTangZhou:2020positive}
	H.-L. Liao, T. Tang and T. Zhou, Positive definiteness of real quadratic forms resulting from variable-step approximations of convolution operators, arXiv:2011.13383v1, 2020.
	
	\bibitem{LiaoZhang:2021}
	H.-L. Liao and Z. Zhang, Analysis of adaptive BDF2 scheme for diffusion equations, Math. Comput., 90 (2021), pp. 1207--1226.
	
	\bibitem{LiaoZhuWang:2022NMTMA}
	H.-L. Liao, X. Zhu and J. Wang, An adaptive L1 time-stepping scheme preserving a compatible energy law for the time-fractional Allen-Cahn equation, Numer. Math. Theor. Meth. Appl., 15(4) (2022), pp. 1128--1146. 
	
	\bibitem{LubichSloanThomee:1996}
	C. Lubich, I. H. Sloan and V. Thom\'{e}e, Nonsmooth data error estimates for approximations of an evolution equation with a positive-type memory term, Math. Comput., 65 (1996), pp. 1--17.
	
	\bibitem{LyuVong:2022JSC}
	P. Lyu and S. Vong, A symmetric fractional-order reduction method for direct nonuniform approximations of semilinear diffusion-wave equations, J. Sci. Comput., 93 (2022), num. 34. 
	
	 \bibitem{Mainardi2010}
	F. Mainardi, Fractional Calculus and Waves in Linear Viscoelasticity, Imperial College Press, London, 2010.
	
	
	\bibitem{McLeanMustapha:2007}
	W. McLean and K. Mustapha, A second-order accurate numerical method for a fractional wave equation, Numer. Math., 105 (2007), pp. 481--510.
	
	\bibitem{McLeanThomee:1996}
	W. McLean and V. Thom\'{e}e, Discretization with variable time steps of an evolution equation with a positive-type memory term, J. Comput. Appl. Math., 69 (1996), pp. 49--69.
	
	
	\bibitem{McLeanThomee:1993}
	W. McLean and V. Thom\'{e}e, Numerical solution of an evolution equation with a positive-type
	memory term, J. Austral. Math. Soc. Ser., 35 (1993), pp. 23--70.
	
	
	\bibitem{MetzlerKlafter:2000}
	R. Metzler and J. Klafter, Accelerating Brownian motion: A fractional dynamics approach to fast diffusion, Europhys. Lett., 51 (2000), pp. 492--498.
	
	\bibitem{MustaphaMustapha:2010}
	K. Mustapha and H. Mustapha, A second-order accurate numerical method for a semilinear integro-differential equation with a weakly singular kernel, IMA J. Numer. Anal., 30 (2010), pp. 555--578.
	
	\bibitem{Mustapha-sinum2020}
	K. Mustapha, An $L1$ approximation for a fractional reaction-diffusion equation, a second-order error analysis over time-graded meshes, SIAM J. Numer. Anal., 58 (2020), pp. 1319--1338.
	
	\bibitem{Mustapha McLean:2013}
	K. Mustapha and W. McLean, Superconvergence of a discontinuous Galerkin method for fractional diffusion and wave equations, SIAM J. Numer. Anal., 51 (2013), pp. 491--515.
	
	\bibitem{MustaphaSchotzau:2014}
	K. Mustapha and D. Sch\"{o}tzau, Well-posedness of {\em hp}-version discontinuous Galerkin methods for fractional diffusion wave equations, IMA J. Numer. Anal., 34 (2014), pp. 1426--1446.
	
	
	\bibitem{Nigmatullin-pb1984}
	R. R. Nigmatullin, To the theoretical explanation of the ``Universal Response", Physica B, 123 (1984), pp. 739--745.
	
	
	\bibitem{QuanTangYang-2020csiam}
	C. Quan, T. Tang and J. Yang, How to define dissipation-preserving energy for time-fractional
	phase-field equations, CSIAM-AM, 1 (2020), pp. 478--490.
		
	\bibitem{QuanTangYang-2020}
	 C. Quan, T. Tang and J. Yang, Numerical energy dissipation for time-fractional phase-field equations, arXiv:2009.06178v1, 2020.
	
	 	
	\bibitem{QuanTangWangYang:2022}
	C. Quan, T. Tang, B. Wang and J. Yang,
	 A decreasing upper bound of energy for time-fractional	phase-field equations, arXiv:2202.12192v1, 2022.
	 
	
	\bibitem{StynesORiordanGracia:2017}
	M. Stynes, E. O'Riordan, and J. L. Gracia, Error analysis of a finite difference method on graded meshes for a time-fractional diffusion equation, SIAM J. Numer. Anal., 55 (2017), pp. 1057--1079.
	
	\bibitem{TangYuZhou:2019}
	T. Tang, H. Yu and T. Zhou, On energy dissipation theory and numerical stability for time-fractional phase-field equations, SIAM J. Sci. Comput., 41 (2019), pp. A3757--A3778.
	
	
	
\end{thebibliography}
\end{document}